\documentclass[11pt,a4paper,reqno]{amsart}
\usepackage{amsmath,amssymb,amsfonts,amsthm}
\usepackage{dsfont}

\usepackage{a4wide}
\usepackage{color}
\usepackage{pdfsync}
\usepackage{graphicx,subfigure}
\usepackage{hyperref}

\newtheorem{theorem}{Theorem}[section]

\newtheorem{corollary}[theorem]{Corollary}
\newtheorem{proposition}[theorem]{Proposition}

\newtheorem{remark}[theorem]{Remark}
\numberwithin{equation}{section}
\makeatletter

\@addtoreset{equation}{section}
\makeatother

\newcommand{\N}{{\mathbb N}}
\newcommand{\Z}{{\mathbb Z}}

\newcommand{\R}{{\mathbb R}}

\newcommand{\eps}{{\varepsilon}}
\def\cal#1{\mathcal{#1}}

\def\K{\cal K}
\def\Tc{\mathcal{T}}
\def\Tca#1{{\Tc^{\eps}_{#1}}}

\def\curl{{\rm curl }\,}
\def\div{{\rm div }\,}
\def\loc{\, {\rm loc}}
\def\supp{{{\rm supp}\;}}

\title[Impermeability through a perforated domain]{Impermeability through a perforated domain for the incompressible 2D Euler equations}

\author[C. Lacave \& N. Masmoudi]{Christophe Lacave \& Nader Masmoudi}

\address[C. Lacave]{Univ Paris Diderot, Sorbonne Paris Cit\'e, Institut de Math\'ematiques de Jussieu-Paris Rive Gauche, 
UMR 7586, CNRS, Sorbonne Universit\'es, UPMC Univ Paris 06, F-75013, Paris, France.}
\email{christophe.lacave@imj-prg.fr}

\address[N. Masmoudi]{Courant Institute, 251 Mercer St., New York, NY 10012, U.S.A.}
\email{masmoudi@cims.nyu.edu}
\date{\today}

\begin{document}
\maketitle

\begin{abstract}
We study the asymptotic behavior of the motion of an ideal incompressible fluid in a perforated domain. The porous medium is composed of inclusions of size $\varepsilon$ separated by distances $d_\varepsilon$ and the fluid fills the exterior.

If the inclusions are distributed on the unit square, the asymptotic behavior depends on the limit of $\frac{d_{\varepsilon}}\varepsilon$ when $\varepsilon$ goes to zero. If $\frac{d_{\varepsilon}}\varepsilon\to \infty$, then the limit motion is not perturbed by the porous medium, namely we recover the Euler solution in the whole space. On the contrary, if $\frac{d_{\varepsilon}}\varepsilon\to 0$, then the fluid cannot penetrate the porous region, namely the limit velocity verifies the Euler equations in the exterior of an impermeable square.

If the inclusions are distributed on the unit segment then the behavior 
depends on the geometry of the inclusion: it is determined by the limit of $\frac{d_{\varepsilon}}{\varepsilon^{2+\frac1\gamma}}$ where $\gamma\in (0,\infty]$ is related to the geometry of the lateral boundaries of the obstacles. If 
$\frac{d_{\varepsilon}}{\varepsilon^{2+\frac1\gamma}} \to \infty$, then the presence of holes is not felt at the limit, whereas an impermeable wall appears if this limit is zero. Therefore, for a distribution in one direction, the critical distance depends on the shape of the inclusions. 
In particular it is equal to $\varepsilon^3$ for balls.
\end{abstract} 

\tableofcontents

\section{Introduction}

The problem studied in this paper is the behavior of the 2D Euler equations in a porous medium. Although very natural, this problem was studied in very 
 few mathematical papers. Most papers have focused on the Laplace equation \cite{CM82}, on Stokes and Navier-Stokes flows \cite{Allaire90a,Mikelic91,Sanchez80,Tartar80}. Recently, more attention was given to the homogenization of other fluid models such as the compressible Navier-Stokes system \cite{Diaz99,Masmoudi02esaim} and the acoustic system \cite{AMBPT13,DAM12}.

Concerning the incompressible Euler system, noticeable exceptions are the works \cite{LionsMasmoudi,MikelicPaoli}, on a weakly nonlinear Euler flow through a regular grid (balls of radius $\varepsilon$, at distance $\varepsilon$ from one another). In this context, a non-trivial limit is obtained as $\varepsilon$ goes to zero.

In the spirit of \cite{Allaire90b} where the limit was studied for holes which are larger or smaller than the critical sizes of the holes, we have obtained in \cite{BLM} a condition such that the ideal incompressible fluid does not feel the presence of the small inclusions. 
Our goal here is to complete the study started in \cite{BLM}.
 We improve the critical distance obtained 
 in \cite{BLM} and we prove that below this
 distance an impermeable boundary appears.

When the inclusions are distributed in one direction, an important difference with the result about the Stokes and Laplace equations is that the critical value depends on the shape of the inclusions and corresponds to very close inclusions (at distance $\varepsilon^3$ for disks distributed on a segment whereas the inclusions are separated by $1/|\ln \eps|$ in \cite{Allaire90b}). Even if physically, it appears natural that the critical distance for a viscous fluid is larger than for an inviscid flow, it was not clear mathematically because the Euler equations is related to a Laplace problem $\Delta \psi =\omega$ where $\omega$ is the vorticity (which is bounded) and $\psi$ the stream function: the velocity verifying $u=\nabla^\perp \psi$. Of course, the main difference lies in the boundary condition.

\subsection{The perforated domain}

We denote the shape of the inclusions by $\K$ and the standard assumption is:
\begin{description}
\item[(H1)] $\K \subset[-1,1]^2$ is a simply-connected compact set of $\R^2$ such that $\partial K\in C^{1,\beta}$ for $\beta >0$ is a Jordan curve.
\end{description}

We consider the case where the inclusions have the same shape:
\begin{equation}\label{domain1}
\K_{i,j}^{\eps}:= z_{i,j}^{\varepsilon} + \tfrac\eps 2 \K,
\end{equation}
where the points $z_{i,j}^{\varepsilon}\in \R^2$ should be fixed in order that the inclusions are disjoints. We assume that the inclusions of size $\varepsilon$ are at least separated by a distance $d_\varepsilon$. Namely, for $i,j\in \Z$ and $\varepsilon>0$, we define 
\begin{equation}\label{domain2}
z_{i,j}^{\varepsilon}:=(\tfrac\eps 2+(i-1)(\eps+d_{\varepsilon}),(j-1)(\eps+d_\eps))=(\tfrac\eps 2,0)+(\eps+d_{\eps})(i-1,j-1).
\end{equation}

\begin{figure}[h!t]
\begin{center}
\includegraphics[height=1.5cm]{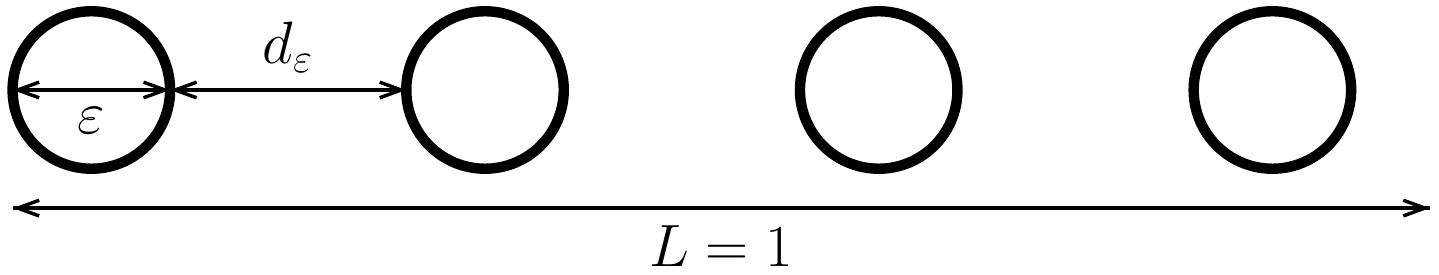}
\caption{Inclusions along the unit segment.}\label{fig.config}
\end{center}
\end{figure}

In the horizontal direction, we consider the maximal number of inclusions that we can distribute on the unit segment $[0,1]$ (see Figure~\ref{fig.config}), hence we consider 
$$i=1,\dots, N_{\varepsilon} \text{ in \eqref{domain1}-\eqref{domain2}},$$
with 
$$N_{\varepsilon}=\left[\frac{1+d_\eps}{\eps+d_\eps}\right]$$
(where $[x]$ denotes the integer part of $x$).

In the vertical direction, we will consider two situations:
\begin{itemize}
 \item inclusions covering the unit square, namely 
 $$j=1,\dots, N_{\varepsilon} \text{ in \eqref{domain1}-\eqref{domain2}};$$
 \item inclusions concentrated on the unit segment, namely 
 $$j=1 \text{ in \eqref{domain1}-\eqref{domain2}}.$$
\end{itemize}

When the obstacles are distributed only in one direction, the critical distance will be much smaller than the size, and the geometry of the lateral boundaries of the inclusion will play a crucial role. We assume that the distance between $\K^\eps_{i,1}$ and $\K^\eps_{i+1,1}$ is reached for $z^\eps_{i,1}+(\tfrac\eps 2,0)$ and $z^\eps_{i+1,1}-(\tfrac\eps 2,0)$. Namely, we assume that $(\pm 1,0)\in \partial \K$, which implies that the distance between two inclusions is $d_{\eps}$:
\[
{\rm d}\Big(\K_{i,1}^{\eps},\K_{i+1,1}^{\eps}\Big)={\rm d}\Big(z_{i,1}^{\varepsilon} + (\tfrac\eps 2,0),z_{i+1,1}^{\varepsilon} - (\tfrac\eps 2,0)\Big)=d_{\varepsilon}.
\]
Of course, we could replace this condition by $(\pm 1,s)\in \partial \K$, for some $s\in [-1,1]$.

We recall that if the boundary 
 $\partial\K$ is locally parametrized around $(1,0)$ by 
\[
x = 1- \rho |y|^{1+\gamma}
\]
with $\rho ,\gamma>0$, then $\gamma$ is called 
the tangency exponent in \cite{CNS,MunnierRamdani}.
 More generally, let us assume that
\begin{description}
\item[(H2)] there exist $\rho_{0}\in (0,1)$, $\rho_{1},\rho_{2}>0$ and $\gamma\in (0,\infty]$ such that
\begin{itemize}
 \item the segment $[-(1-\rho_{1} |s|^{1+\gamma}),(1-\rho_{1} |s|^{1+\gamma})]\times\{s\}\subset \K$ for all $s\in [-\rho_{0},\rho_{0}]$;
 \item the set $\K\cap(\R\times \{s\}) \subset [-(1-\rho_{2} |s|^{1+\gamma}),(1-\rho_{2} |s|^{1+\gamma})]\times\{s\}$ for all $s\in [-1,1]$.
\end{itemize}
\end{description}
In particular, (H2) encodes the fact that $(\pm 1,0)\in \partial \K$, which will be used when the inclusions are distributed on the segment. For inclusions covering the square, our main result is independent of $\rho_{0}$, $\rho_{1}$ and $\gamma$, but we will need $(\pm 1,0), (0,\pm 1) \in \partial \K$ (see Section~\ref{sect shape} for possible extensions).

For example, if the inclusion is the unit ball, $\partial \K$ is locally parametrized around $(1,0)$ by $(\sqrt{1-s^2},s)\sim_{0} (1-\frac12 s^2,s)$, i.e. we can choose in this case $\gamma =1$ (with $\rho_{2}<1/2<\rho_{1}$). Note that $\gamma=\infty$ corresponds to the case where the solid is flat near $(\pm 1,0)$: $[( 1,-\rho_{0}),( 1,\rho_{0})]\subset \partial \K$.

Let us note that $\gamma$ does not correspond to the 
radius of curvature $R_{c}$: for $\gamma \in (0,1)$ the radius $R_{c}=0$, for $\gamma \in (1,\infty]$ we have $R_{c}=\infty $, and for $\gamma=1$ we compute $R_{c}=1/(2\rho)$.

Throughout the paper, the fluid domain will be the exterior of the inclusions:
\begin{equation}\label{Omega2eps}
 \Omega^{\varepsilon}_{2}:=\R^2 \setminus \Big( \bigcup_{i=1}^{N_{\eps}} \bigcup_{j=1}^{N_{\eps}} \K_{i,j}^{\varepsilon}\Big),
\end{equation}
or
\begin{equation}\label{Omega1eps}
\Omega^{\varepsilon}_{1}:=\R^2 \setminus \Big( \bigcup_{i=1}^{N_{\eps}} \K_{i,1}^{\varepsilon}\Big).
\end{equation}

\subsection{The Euler equations and former results}

The velocity $u^\eps:=(u^\eps_{1},u^\eps_{2})(t,x)$ of an ideal incompressible fluid filling a domain $\Omega^\eps$ is governed by the Euler equations
\begin{eqnarray}
\partial_{t} u^\eps + u^\eps\cdot \nabla u^\eps = -\nabla p^\eps, && (t,x)\in (0,\infty)\times \Omega^\eps; \label{Euler1}\\
\div u^\eps = 0, && (t,x)\in [0,\infty)\times \Omega^\eps;\label{Euler2}\\ 
u^\eps\cdot n = 0, && (t,x)\in [0,\infty)\times \partial \Omega^\eps;\label{Euler3}\\
u^\eps(0,\cdot)=u^\eps_{0}, && x\in \Omega^\eps,\label{Euler4}
\end{eqnarray}
where $p^\eps$ is the pressure. In dimension two, the natural quantity for these equations is the vorticity:
\[
\omega^\eps:= \curl u^\eps = \partial_{1} u^\eps_{2}-\partial_{2}u^\eps_{1}
\]
because taking formally the curl of \eqref{Euler1}, we note that $\omega^\eps$ verifies a transport equation:
\begin{equation}
\partial_{t} \omega^\eps + u^\eps\cdot \nabla \omega^\eps = 0, \quad (t,x)\in (0,\infty)\times \Omega^\eps . \label{Euler5}
\end{equation}

In smooth domains, the well-posedness of the Euler equations was established long time ago: the existence and uniqueness of global strong solution for smooth initial data is a result of Wolibner \cite{Wolibner}, whereas the existence and uniqueness of global weak solutions for initial data with bounded vorticity was
 obtained by Yudovich \cite{Yudo} (see Kikuchi \cite{Kikuchi} for exterior domains).

To fix an initial data for all $\varepsilon$, the standard way is to fix an initial vorticity $\omega_{0}\in L^\infty_{c}(\R^2)$ (i.e. $\omega_{0}$ bounded and compactly supported) and to consider the unique initial velocity $u_{0}^\varepsilon$ such that:
\begin{equation}\label{initial}
\begin{split}
\div u_{0}^\varepsilon =0 \text{ in } \Omega^{\varepsilon},\quad \curl u_{0}^\varepsilon =\omega_{0}\text{ in } \Omega^{\varepsilon}, \quad u_{0}^\varepsilon \cdot n = 0 \text{ on } \partial \Omega^{\varepsilon}, \\
\lim_{|x|\to\infty}u_{0}^\eps(x)=0,\quad \oint_{\partial \K_{i,j}^{\varepsilon}} u_{0}^\varepsilon\cdot \tau\, ds=0 \text{ for all }i,j.
\end{split}
\end{equation}

{\it The goal of this article is to determine the limit of the Euler solutions when $\varepsilon$ and $d_{\varepsilon}$ tend to zero. 
More precisely, depending on the distribution of the 
inclusions, we prove the convergence to the solution of Euler in the whole plane, in the exterior of the unit segment or in the exterior of the unit square.}

The asymptotic behavior of the Euler equations around inclusions which shrink to points was the subject of several recent works: the case of one shrinking obstacle in the whole plane \cite{ILL}, of one shrinking obstacle in a bounded domain \cite{Lopes}, of an infinite number of shrinking disks \cite{LLL}. In all these works, the authors have also treated the case where the initial circulations are non zero, which adds some difficulties and provides interesting asymptotics. 
Nevertheless, in all these cases the distance between the obstacles is large, and if the initial circulations are zero, 
the theorems therein can be understood by stating that the limit motion is not perturbed by the shrinking inclusions.

\subsection{Main result}

In \cite{BLM} we have looked for the distances $d_{\eps}= \varepsilon^\alpha$ between obstacles 
($\alpha>0$) for which we do not feel the presence of the inclusions. We have proved therein that for $\alpha<1$ when inclusions cover the square and for $\alpha<2$ when inclusions are distributed in one direction, then the solutions $(u^\eps ,\omega^\eps)$ of the Euler equations on $\Omega^{\varepsilon}$ with initial velocity $u_{0}^\eps$ \eqref{initial} converge to the unique global solution to the Euler equations in the whole plane $\R^2$, with initial vorticity $\omega_{0}$.

When inclusions cover the square, the first goal of this paper is to observe, for $\alpha>1$, the presence of an impermeable square at the limit, i.e. we get $u\cdot n=0$ on the boundary of the square, no matter what the geometry of the inclusions is.

When inclusions are distributed on the segment, we will prove in the second theorem that the critical value of $\alpha$ depends on the geometry:
\[
\alpha_{c}(\gamma)={2+\frac1\gamma}.
\]
For instance:
\begin{itemize} 
 \item if inclusions are asymptotically vertically flat ($\gamma=\infty$ in (H2)), an impermeable segment appears for $\alpha>\alpha_{c}(\infty)=2$;
\item for circular inclusions, we prove the impermeability of the porous medium for $\alpha>\alpha_{c}(1)=3$ and prove permeability for $\alpha<\alpha_{c}(1)=3$ improving 
 the analysis performed in \cite{BLM}.
\end{itemize} 
 
A main novelty of this work is to show that the critical distance depends on the geometry of the inclusions (in the case of a segment),
 which is not the case for the Laplace, Stokes and Navier-Stokes problems.

Another slight improvement compared to \cite{BLM}, included here, is to consider more general $d_{\eps}$ than power of $\varepsilon$. This requires us to be more precise and not even 
losing some $log$ in some of the estimates in Section~\ref{sec:permeability}.

\subsubsection{Case of the square}
When the inclusions cover a square, the criterion is independent of the geometry: 

\begin{theorem} \label{main2} Assume that $\K$ verifies {\rm (H1)}. Let $\omega_{0}\in L^\infty_{c}(\R^2)$ and $(u^\eps,\omega^\eps)$ be the global weak solution to the Euler equations \eqref{Euler1}-\eqref{Euler5} on 
$$\Omega^{\varepsilon}_{2}=\R^2 \setminus \Big( \bigcup_{i=1}^{N_{\eps}} \bigcup_{j=1}^{N_{\eps}} \K_{i,j}^{\varepsilon}\Big)\quad \text{(with $ \K_{i,j}^{\varepsilon}$ defined in \eqref{domain1}-\eqref{domain2})},$$
with initial vorticity $\omega_{0}\vert_{\Omega^{\varepsilon}}$ and initial circulations $0$ around the inclusions (see \eqref{initial}).
\begin{enumerate}
 \item[(i)] If 
 $$
 \frac{d_{\eps}}{\eps} \to \infty \quad \text{ for a sequence }\varepsilon \to 0
 $$
 then
 \begin{itemize}
\item $u^\eps \to u$ strongly in $L^2_{\loc}(\R^+\times\R^2)$ and $\omega^\eps \rightharpoonup {\omega}$ weak $*$ in $L^\infty(\R^+\times\R^2)$;
\item the limit pair $(u,{\omega})$ is the unique global solution to the Euler equations in the whole plane $\R^2$, with initial vorticity $\omega_{0}$.
\end{itemize}
 \item[(ii)] If 
 $$
 \frac{d_{\eps}}{\eps} \to 0 \quad \text{ for a sequence }\varepsilon \to 0
 $$
 and that $(\pm 1,0),(0,\pm 1)\in \partial \K$, then there exists a subsequence such that
 \begin{itemize}
\item $u^{\eps} \rightharpoonup u$ weak $*$ in $L^\infty_{\loc}(\R^+; L^2_{\loc}(\R^2\setminus [0,1]^2))$ and $ \omega^{\eps} \rightharpoonup {\omega}$ weak $*$ in $L^\infty(\R^+\times(\R^2 \setminus [0,1]^2))$;
\item the limit pair $(u,{\omega})$ is a global weak solution to the Euler equations in $\R^2\setminus [0,1]^2$, with $u\cdot n=0$ on the boundary, with initial vorticity $\omega_{0}$ and initial circulation $0$ around the square.
\end{itemize}
\end{enumerate}
\end{theorem}
The tangency condition is verified where $n$ is defined, i.e. on the boundary of the square except the corners: for all $x\in \Big((0,1)\times \{0,1\}\Big)\cup \Big( \{0,1\}\times(0,1)\Big)$. Let us note that it does not imply that the velocity vanishes at the corners, because the velocity may blow up near the corners (see Remark~\ref{rem Grisvard} for more details).

The situation
\[
 \frac{d_{\eps}}{\eps} \to C >0,
\]
is related to the problem considered in \cite{LionsMasmoudi,MikelicPaoli}. In these works, the authors considered the case where the inclusions cover a bounded domain $\widetilde \Omega$ with $d_{\eps}=\eps$. As there are inclusions everywhere in the fluid domain, they have assumed that the velocity is small $u^\eps=\eps \tilde u^\eps$ and studied the limit when $\tilde u^\eps$ verifies a weakly non-linear equation
\[
\partial_{t} \tilde u^\eps +\eps \tilde u^\eps\cdot \nabla \tilde u^\eps = -\nabla \tilde p^\eps.
\]
Using the notion of two-scale convergence, they obtained that the limit of $\tilde u^\eps$ verifies an homogenized system which couples a cell problem with the macroscopic one.

\subsubsection{Case of the segment}
When the inclusions are distributed only in one direction, the critical distance is much smaller than the size, and depends on the exponent $\gamma$ of the lateral boundaries of the obstacle. Namely, the result reads as follows 
\begin{theorem} \label{main1} Assume that $\K$ verifies {\rm (H1)-(H2)} with $\gamma\in (0,\infty]$. 
Let $\omega_{0}\in L^\infty_{c}(\R^2)$ and $(u^\eps,\omega^\eps)$ be the global weak solution to the Euler equations \eqref{Euler1}-\eqref{Euler5} on 
$$\Omega^{\varepsilon}_{1}=\R^2 \setminus \Big( \bigcup_{i=1}^{N_{\eps}} \K_{i,1}^{\varepsilon}\Big)\quad \text{(with $ \K_{i,j}^{\varepsilon}$ defined in \eqref{domain1}-\eqref{domain2})},$$
with initial vorticity $\omega_{0}\vert_{\Omega^{\varepsilon}}$ and initial circulations $0$ around the inclusions (see \eqref{initial}).
\begin{enumerate}
 \item[(i)] If 
 $$
 \frac{d_{\eps}}{\eps^{2+\frac1\gamma}} \to \infty \quad \text{ for a sequence }\varepsilon \to 0
 $$
 then
 \begin{itemize}
\item $u^\eps \to u$ strongly in $L^2_{\loc}(\R^+\times\R^2)$ and $\omega^\eps \rightharpoonup {\omega}$ weak $*$ in $L^\infty(\R^+\times\R^2)$;
\item the limit pair $(u,{\omega})$ is the unique global solution to the Euler equations in the whole plane $\R^2$, with initial vorticity $\omega_{0}$.
\end{itemize}
 \item[(ii)] If 
 $$
 \frac{d_{\eps}}{\eps^{2+\frac1\gamma}} \to 0 \quad \text{ for a sequence }\varepsilon \to 0
 $$
 then there exists a subsequence such that
 \begin{itemize}
\item $u^{\eps} \rightharpoonup u$ weak $*$ in $L^\infty_{\loc}(\R^+; L^2_{\loc}(\R^2\setminus ([0,1]\times\{0\})))$ and $ \omega^{\eps} \rightharpoonup {\omega}$ weak $*$ in $L^\infty(\R^+\times(\R^2 \setminus ([0,1]\times\{0\})))$;
\item the limit pair $(u,{\omega})$ is a global weak solution to the Euler equations in $\R^2\setminus ([0,1]\times\{0\})$, with $u\cdot n=0$ on the boundary, with initial vorticity $\omega_{0}$ and initial circulation $0$ around the segment.
\end{itemize}
\end{enumerate}
\end{theorem}

The tangency condition is verified on both sides of the segment, except at the endpoints where the velocity may blow up (see Remark~\ref{rem Grisvard}). As it will be explained just below Remark~\ref{rem Grisvard}, the explicit behavior of $u_{0}$ gives the integrability of the trace and the circulation $\oint_{\partial (\R^2\setminus ([0,1]\times\{0\}))} u_{0}\cdot \tau \, ds$ should be understood as the integral on the upper value $\int_{0}^1 u_{0}^{up} \cdot (-e_{1}) $ plus the lower value $\int_{0}^1 u_{0}^{down} \cdot e_{1}$. 

\begin{remark}
The extraction of a subsequence appears only in Section~\ref{sec2}, due to general compactness results (namely, we apply the Banach-Alaoglu's theorem to the vorticity). 

In the case (i) of Theorems~\ref{main2} and \ref{main1}, we recover at the limit the solution of the Euler equations in the whole plane. By the Yudovich's theorem, this solution is unique and we do not need to extract a subsequence.

Unfortunately, the uniqueness result is a hard issue in the exterior of the square or the segment, because in such domains, the velocity does not belong to $\cap_{p\in [2,\infty)} W^{1,p}$. Actually, $u$ is not even bounded and blows up in the vicinity of the corners (like $1/\sqrt{|x|}$ close to the endpoints of the segment and like $1/|x|^{1/3}$ close to the corners of a square). We refer to \cite{Lac-uni} for partial results in this setting.

The $weak\ *$ limit for $u^\varepsilon$ in the case (ii) is not enough to pass to the limit in the non-linear term $(u^\eps\otimes u^\eps):\nabla \psi$ (or $\omega^\eps u^\eps\cdot \nabla \varphi$ in the vorticity formulation), and we will establish later (see Step 4 in Section~\ref{sec-cvg-small}) that for any $T>0$ and $\mathcal{O}\Subset \R^2 \setminus [0,1]^2$ (or $\mathcal{O}\Subset \R^2 \setminus ([0,1]\times \{ 0\})$), we have
\[
\mathbb{P}_{\mathcal{O}} u^\varepsilon \to \mathbb{P}_{\mathcal{O}} u \text{ strongly in } L^\infty(0,T;L^2(\mathcal{O})),
\]
where $\mathbb{P}_{\mathcal{O}}$ is the Leray projector.
\end{remark}

The main novelty of this theorem is to get a criterion which depends on the geometry of the inclusions. This is natural because the distance is much smaller than the size of the holes.
 
For the Laplace, Stokes and Navier-Stokes problem in dimension two with Dirichlet boundary condition, the critical criterion separating between impermeability and permeability 
 is whether $d_{\eps} \sqrt{|\ln \eps|}$ goes to $\infty$ or $0$ when inclusions cover a bounded domain $\widetilde \Omega$ and whether $d_{\eps} |\ln \eps|$ goes to $\infty$ or $0$ when they are distributed in one direction. For these systems, the main argument is an energy argument. Roughly, to prove the convergence of $u^\eps$ to the solution $u$ of the Laplace problem in the whole domain, one introduces $W^\eps=u^\eps - \varphi^\eps u$ where $\varphi^\eps$ is a cutoff function vanishing at the boundaries of the inclusions. To pass to the limit in the weak formulation, the main difficulty is to study terms like $\int \nabla( \varphi^\eps u) : \nabla( \varphi^\eps \psi)$ where $\psi$ is a test function. Choosing the cutoff function which minimize the $L^2$ norm of the gradient (see Remark~\ref{rem:cutoff}), we find
\[
\| \nabla \varphi^\eps \|_{L^2} = \frac C{d_{\eps}^{\mu} \sqrt{|\ln \eps |}}
\]
with $\mu=1$ if the inclusions are distributed on a surface and $\mu=1/2$ for a distribution on a curve. Therefore, for inclusions covering a bounded domain $\widetilde \Omega$ and when $d_{\eps} \sqrt{|\ln \eps|}\to C>0$, we obtain at the limit an homogenized system: see Cioranescu-Murat \cite{CM82} for the Laplace equation and Allaire \cite{Allaire90a} for the Stokes and Navier-Stokes systems. If $d_{\eps} \sqrt{|\ln \eps|}\to \infty$, Allaire has proved in \cite{Allaire90b} that the limit is not perturbed by the porous medium, and he has also proved the homogenized phenomenon (a Brinkmann-type law) when $d_{\eps} |\ln \eps|\to C>0$ for inclusions distributed on a curve.

Let us note that the case $d_{\eps} |\ln \eps|\to 0$ for inclusions distributed on a curve is not treated by Allaire \cite{Allaire90b}. With smaller distance than the critical value, he only considered the case of inclusions covering a bounded domain $\widetilde \Omega$, and as there are inclusions everywhere, the velocity is assumed small and he gets a limit for the quantity $u^\eps/d_{\eps} \sqrt{|\ln \eps|}$ in order to pass to the limit in the term of the type
\[
\int \nabla\left(\frac{u^\eps}{d_{\eps} \sqrt{|\ln \eps|}}\right): (d_{\eps} \sqrt{|\ln \eps|} \nabla \varphi^\eps \otimes \psi).
\]
In this setting, he obtained Darcy's law at the limit. Actually, Allaire needs to assume that the size of the holes is still smaller than the inter-hole distance (see \cite[Equation (3.4.2)]{Allaire90b}), hence the case $\frac{d_{\varepsilon}}\varepsilon\to 0$ in our main results is also a novelty of the present paper and \cite{BLM}.

For inclusions distributed on a curve, there are also many works in the periodic setting $d_{\eps}=\eps$ (through grids, sieves or porous walls) we refer e.g. to \cite{conca1, conca2,SP1}.

The reminder of this paper is composed of four parts. In the following section, we reduce our problem to a time-independent problem. First, we derive some a priori estimates which allow us to pass to the limit far away from the porous medium. Next, we write two key propositions dealing with the behavior of divergence free vector fields in a porous medium. We will prove that these propositions imply Theorems~\ref{main2} and \ref{main1}.

The core of our analysis is the proof of these two propositions where the quantities $\frac{d_{\eps}}\eps$ (for the square) and $\frac{d_{\eps}}{\eps^{2+\frac1\gamma}}$ (for the segment) 
 will appear naturally as critical quantities to get permeability vs impermeability.

In Section~\ref{sec:permeability}, we prove Proposition~\ref{prop:key1} concerning permeability for
 inclusions separated by large enough distances. We will follow therein the strategy in \cite{BLM}, which is somehow related to homogenization ideas. The cell problem is the exterior of one inclusion where the velocity verifies the tangency condition. We introduce and study a correction of the solution in the whole plane. This correction is based on the explicit Green kernel formula in the exterior of one simply connected compact set, via Riemann mappings. To localize our formula, we need to truncate the stream function by introducing some cutoff functions $\varphi^\eps$ and we will discuss the optimal choice of $\varphi^\eps$.

Section~\ref{sec:impermeability} deals with Proposition~\ref{prop:key2} concerning
 impermeability for inclusions separated by small enough distances.
To prove that a wall appears, the main point is to get the impermeability condition on the boundary of the porous medium. The approach is totally different from the 
permeability case. Here, we use an $L^2$ argument and the computation of the area between two obstacles, which depends on the coefficient $\gamma$. Indeed, for small distance, it makes sense, physically, that the fluid crosses easier between circular inclusions ($\gamma=1$) than between squares ($\gamma=\infty$). 

Finally, Section~\ref{sect final} is dedicated to further discussions. We will give some extensions on the inclusion geometry where our results hold. For instance, note that $\gamma =0$ in (H2) would correspond to the case where $\partial \K$ is a corner of angle $2{\rm Arctan}\ \frac1{\rho_{1}}$, which would not verify $\partial \K\in C^1$. This case is not treated in the main theorems, and we refer to Section~\ref{sect shape} for a discussion about inclusions with corners. We do not treat the case of a cusp, where we could only take $\gamma \in (-1,0)$. 
Section~\ref{sect final} contains also a consequence on the Leray projection: we deduce from our analysis that the Leray projection is not uniformly continuous from $L^p$ to $L^p$ (with $p\in [1,2)$) in porous medium where $d_{\eps} \ll \eps$.

\medskip

If it appears natural that the critical distance is smaller for ideal fluid than for viscous fluid, this is not so clear mathematically. We could argue that the Euler equations and the Laplace, Stokes or Navier-Stokes problems
are not at same order, but in practical, the velocity comes from a Laplace problem $\Delta \psi^\eps = \omega_0$ where $\psi^\eps$ is the stream function, and we are looking for $L^2$ estimates of $\nabla^\perp \psi^\eps = u^\eps_{0}$. Then, we could think that the difference comes from the boundary condition and that the tangency and the zero circulation conditions are less restricting for ideal fluid than the Dirichlet boundary condition for viscous fluid. Nevertheless, this is not a sufficient justification: Allaire showed in \cite{Allaire91} that the scaling is the same for Navier-Stokes equations for the two types of boundary conditions (Dirichlet/Robin). We finish this introduction by referring to a recent result of Lacave and Mazzucato \cite{LaMa}, where one considers the limit of the porous medium together with the vanishing viscosity problem.

{\it Notations: } throughout the paper, $\Omega_{2}^\eps$ denotes the exterior domain when inclusions are distributed in the two directions \eqref{Omega2eps} whereas $\Omega_{1}^\eps$ concerns the case of inclusions distributed in one direction \eqref{Omega1eps}. We will use the notation $\Omega^\eps$ when an argument holds for both $\Omega_{1}^\eps$ and $\Omega_{2}^\eps$. For shortness, we also denote the exterior of the unit square and of the unit segment by
\[
\Omega_{2}:=\R^2 \setminus [0,1]^2 \quad \text{and} \quad \Omega_{1}:= \R^2 \setminus ([0,1]\times \{ 0\}),
\]
and $\Omega$ for arguments holding for both.

\section{Proof of convergence}\label{sec2}

In this section, we derive some standard a priori estimates for solutions of the Euler equations. Next, we give the two key propositions and we prove that Theorems~\ref{main2} and \ref{main1} will follow.

\subsection{A priori estimates for Euler solutions}\label{subsect initial}

Let $\omega_{0}\in L^\infty_{c}(\R^2)$ be given. By standard results related to the Hodge-De Rham theorem, there exists a unique vector fields $u_{0}^\eps\in L^2_{\loc}(\overline{\Omega^{\varepsilon}})$ verifying \eqref{initial} (see e.g. \cite{Kikuchi} for a proof in smooth domain and \cite{GV-Lac2} for any irregular domain).

Let us denote
\[
\omega_{0}^\eps := \omega_{0} \mathds{1}_{\Omega^\eps}
\]
then the vector field defined in the whole plane
\[
K_{\R^2}[\omega_{0}^\eps](x):=\frac1{2\pi}\int_{\R^2}\frac{(x-y)^\perp}{|x-y|^2} \omega_{0}^\eps(y)\, dy
\]
is the unique solution of 
\[
\div K_{\R^2}[\omega_{0}^\eps]=0,\quad \curl K_{\R^2}[\omega_{0}^\eps]=\omega_{0}^\eps,\quad \lim_{|x|\to \infty }K_{\R^2}[\omega_{0}^\eps](x)=0.
\]
We choose $R$ so that $\partial \Omega^\eps \subset B(0,R)$, and $\chi$ a smooth cutoff function such that $\chi \equiv 0$ in $B(0,R)$ and $\chi \equiv 1$ in $B(0,R+1)^c$. On the one hand, as $\omega_{0}$ is compactly supported, we have
\[
\widehat{u}_{0}^\eps:= K_{\R^2}[\omega_{0}^\eps] - \frac{\int \omega_{0}^\eps}{2\pi} \nabla^\perp ( \chi \ln |x|) \text{ is uniformly bounded in } L^2(\R^2).
\]
This vector field satisfies:
\begin{gather*}
\div \widehat{u}_{0}^\eps =0 \text{ in } \Omega^{\varepsilon},\quad \curl \widehat{u}_{0}^\eps =\omega_{0} - \frac{\int \omega_{0}^\eps}{2\pi} \Delta ( \chi \ln |x|) \text{ in } \Omega^{\varepsilon}, \\
\lim_{|x|\to\infty}\widehat{u}_{0}^\eps(x)=0,\quad \oint_{\partial \K_{i,j}^{\varepsilon}} \widehat{u}_{0}^\eps\cdot \tau\, ds=0 \text{ for all }i,j,
\end{gather*}
where we have used that $\oint_{\partial \K_{i,j}^{\varepsilon}} K_{\R^2}[\omega_{0}^\eps]\cdot \tau\, ds = \int_{\K_{i,j}^{\varepsilon}}\omega_{0}^\eps=0$ (Stokes formula).

On the other hand, the vector field
\[
\widetilde{u}_{0}^\eps:= u_{0}^\eps- \frac{\int \omega_{0}^\eps}{2\pi} \nabla^\perp ( \chi \ln |x|) 
\]
verifies
\begin{gather*}
\div \widetilde{u}_{0}^\eps =0 \text{ in } \Omega^{\varepsilon},\quad \curl \widetilde{u}_{0}^\eps =\omega_{0} - \frac{\int \omega_{0}^\eps}{2\pi} \Delta ( \chi \ln |x|) \text{ in } \Omega^{\varepsilon}, \quad \widetilde{u}_{0}^\varepsilon \cdot n = 0 \text{ on } \partial \Omega^{\varepsilon},\\
\lim_{|x|\to\infty}\widetilde{u}_{0}^\eps(x)=0,\quad \oint_{\partial \K_{i,j}^{\varepsilon}} \widetilde{u}_{0}^\eps\cdot \tau\, ds=0 \text{ for all }i,j,
\end{gather*}
which implies that $\widetilde{u}_{0}^\eps$ is the Leray projection\footnote{projection on divergence free vector fields which are tangent to the boundary.} of $\widehat{u}_{0}^\eps$. By orthogonality of this projection in $L^2$, we get that
\begin{equation}\label{est-utilde0}
 \| \widetilde{u}_{0}^\eps \|_{L^2(\Omega^\eps)} \leq \| \widehat{u}_{0}^\eps \|_{L^2(\R^2)} \text{ is uniformly bounded.}
\end{equation}

For every $\eps$, it is well-known that there exists a global weak solution 
\[
u^\eps \in L^\infty_{\loc}(\R^+;L^2_{\loc}(\overline{\Omega^\eps})) \quad \text{and} \quad \omega^\eps \in L^\infty(\R^+;L^1\cap L^\infty(\Omega^\eps))
\]
to the Euler equations \eqref{Euler1}-\eqref{Euler5} (see \cite{Kikuchi} for smooth domain and \cite{GV-Lac} for irregular domain). The weak formulation of \eqref{Euler1} and \eqref{Euler4} reads:
\begin{equation}\label{weak Euler}
\int_{0}^\infty\int_{\Omega^\eps} u^\varepsilon \cdot \partial_{t} \phi\, dxdt+\int_{0}^\infty\int_{\Omega} (u^\varepsilon \otimes u^\varepsilon):\nabla \phi\, dxdt=-\int_{\Omega^\eps} u_{0}^\varepsilon\cdot \phi(0,\cdot) \, dx,
\end{equation}
for all divergence free test function $\phi\in C^\infty_{c}([0,\infty)\times \Omega^\eps)$.
This solution verifies 
\begin{itemize}
 \item the transport equation \eqref{Euler5} in the sense of distribution;
 \item the inequality
\begin{equation}\label{om-cons}
 \| \omega^\eps \|_{L^p(\Omega^\eps)} \leq \| \omega^\eps_{0} \|_{L^p(\Omega^\eps)}\leq \| \omega_{0} \|_{L^p(\R^2)}, \quad \forall p\in [1,\infty];
\end{equation}
\item the conservation of the circulation around $\K^\eps_{i,j}$ for all $i,j$ (the Kelvin theorem).
\end{itemize}
Moreover, if the domain is smooth, this solution is unique and is a renormalized solution of \eqref{Euler5} in the sense of DiPerna-Lions.

As $(u^\eps,\omega^\eps)$ is defined on an $\varepsilon$-dependent domain, we extend these functions by zero in the holes. All the convergence results in this paper concern the extended functions. 

Now, let us derive an estimate of $u^\eps$ in $L^\infty_{\loc}(\R^+;L^2_{\loc}(\overline{\Omega}))$ uniformly in $\varepsilon$. The function
\[
\widetilde{u}^\eps:= u^\eps- \frac{\int \omega_{0}^\eps}{2\pi} \nabla^\perp ( \chi \ln |x|)
\]
is divergence free, tangent to the boundary and verifies in a weak sense
\[
\partial_{t} \widetilde{u}^\eps + u^\eps \cdot \nabla \widetilde{u}^\eps = -\nabla p^\eps - \Big(\widetilde{u}^\eps + \frac{\int \omega_{0}^\eps}{2\pi} \nabla^\perp ( \chi \ln |x|) \Big) \cdot\nabla \Big( \frac{\int \omega_{0}^\eps}{2\pi} \nabla^\perp ( \chi \ln |x|)\Big), \quad \widetilde{u}^\eps(0,\cdot)= \widetilde{u}^\eps_{0}\in L^2(\Omega^\eps).
\]
Multiplying by $\widetilde{u}^\eps$ and integrating (or considering a test function $\phi = \text{\Large $\mathds{1}$}_{[0,t]} \widetilde{u}^\eps$, which is possible by density), we obtain
\begin{align*}
 \| \widetilde{u}^\eps (t,\cdot) \|_{L^2(\Omega^\eps)}^2 & \leq \| \widetilde{u}^\eps_{0} \|_{L^2(\Omega^\eps)}^2 + C \int_{0}^t\Big( \| \widetilde{u}^\eps(s,\cdot) \|_{L^2(\Omega^\eps)}^2 +\| \widetilde{u}^\eps(s,\cdot) \|_{L^2(\Omega^\eps)}\Big)\, ds\\
 & \leq \| \widetilde{u}^\eps_{0} \|_{L^2(\Omega^\eps)}^2 + C \int_{0}^t \Big(\| \widetilde{u}^\eps(s,\cdot) \|_{L^2(\Omega^\eps)}^2+1\Big)\, ds
\end{align*}
and we deduce by Gronwall inequality:
\[
 \| \widetilde{u}^\eps (t,\cdot) \|_{L^2(\Omega^\eps)}^2 \leq (\| \widetilde{u}^\eps_{0} \|_{L^2(\Omega^\eps)}^2+1) e^{Ct} .
\]
Putting together with the uniform estimate of $\widetilde{u}^\eps_{0}$ \eqref{est-utilde0}, we get that $\widetilde{u}^\eps$ is uniformly bounded in $L^\infty_{\loc}(\R^+;L^2(\Omega))$, hence
\begin{equation}\label{estL2:ueps}
 u^\eps \text{ is uniformly bounded in }L^\infty_{\loc}(\R^+;L^2_{\loc}(\overline{\Omega})).
\end{equation}

\begin{remark}
The assumption concerning the zero circulation of $u_{0}^\varepsilon$ is crucial to obtain a uniform $L^2$ estimate. For non-zero circulations $\gamma_i$, the authors in \cite{LLL} have had to develop an $L^p$ theory for $p<2$: when an inclusion shrinks to a point $x_{i,j}$, a reminiscent term appears of the form $\gamma_i \frac{(x-x_{i,j})^\perp}{|x-x_{i,j}|^2}$, which belongs to $L^p_{\loc}$ only for $p<2$ (see also \cite{ILL}). However, the continuity of the Leray projector in $L^p$ ($p<2$) uniformly in $\varepsilon$ is unclear (see Section~\ref{sect final} for further discussions).
\end{remark}

\begin{remark}
 As the domain is regular (namely the boundary is $C^{1,\beta}$), elliptic estimates imply that $u^\eps$ is continuous up to the boundary, which gives a meaning (in a strong sense) to the circulation and tangency conditions.
\end{remark}

\subsection{Convergence for large distance}
 
The analysis in \cite{BLM} is based on the remark that $u^\eps$ satisfies at each time
 \begin{gather*}
\div u^\eps =0 \text{ in } \Omega^{\varepsilon},\quad \curl u^\eps =\omega^\eps \text{ in } \Omega^{\varepsilon}, \quad u^\eps\cdot n =0 \text{ on } \partial\Omega^{\varepsilon}\\
\lim_{|x|\to\infty}u^\eps(t,x)=0,\quad \oint_{\partial \K_{i,j}^{\varepsilon}} u^\eps\cdot \tau\, ds=0 \text{ for all }i,j,
\end{gather*}
whereas $K_{\R^2}[\omega^\eps]$ verifies
 \begin{gather*}
\div K_{\R^2}[\omega^\eps] =0 \text{ in } \Omega^{\varepsilon},\quad \curl K_{\R^2}[\omega^\eps] =\omega^\eps \text{ in } \Omega^{\varepsilon}, \\
\lim_{|x|\to\infty}K_{\R^2}[\omega^\eps](t,x)=0,\quad \oint_{\partial \K_{i,j}^{\varepsilon}} K_{\R^2}[\omega^\eps]\cdot \tau\, ds=0 \text{ for all }i,j.
\end{gather*}
The main idea is to introduce, for any function $f\in L^\infty_{c}(\overline{\Omega^\eps})$, an explicit 
approximate solution $v^\eps[f]$ such that
\begin{equation}\label{eq.veps}
\div v^\eps[f] =0 \text{ in } \Omega^{\varepsilon}, \quad v^\eps[f]\cdot n =0 \text{ on }\partial \Omega^{\varepsilon}, \quad \lim_{|x|\to\infty}v^\eps[f](t,x)=0
\end{equation}
which is close to $K_{\R^2}[f]$ in $L^2$. More precisely, the core of the permeability result in Theorems~\ref{main2} and \ref{main1} is contained in the following key proposition: 
 
\begin{proposition} \label{prop:key1} {\bf (Permeability)} 
Assume that $\K$ verifies {\rm (H1)}. 
\begin{itemize} 
 \item Consider $\Omega_{2}^\eps$ \eqref{Omega2eps} in the case $d_{\eps} \geq \eps$. For any $f\in L^\infty_{c}(\overline{\Omega_{2}^\eps})$ there exists $v^\eps[f]$ satisfying \eqref{eq.veps} such that
\begin{equation*}
\|K_{\R^2}[f] - v^\eps [f] \|_{L^2(\Omega^\eps_{2})}\leq C \| f \|_{L^1\cap L^\infty} \frac{\eps} {d_{\eps}},
\end{equation*}
with $C$ independent of $f$ and $\eps$.
 \item Consider $\Omega_{1}^\eps$ \eqref{Omega1eps} with $\K$ verifying also {\rm (H2)}. For any $f\in L^\infty_{c}(\overline{\Omega_{1}^\eps})$ there exists $v^\eps[f]$ satisfying \eqref{eq.veps} such that
\begin{equation*}
\|K_{\R^2}[f] - v^\eps [f] \|_{L^2(\Omega^\eps_{1})}\leq C \| f \|_{L^1\cap L^\infty} \sqrt{\eps} \Big(1 + \Big(\frac{\eps}{d_{\eps}}\Big)^{\frac{\gamma}{2(\gamma+1)}}\Big),
\end{equation*}
with $C$ independent of $f$ and $\eps$.
\end{itemize}
\end{proposition}

This proposition is proved in Section~\ref{sec:permeability}. The main idea 
is to estimate 
 the correction of the tangency condition in the exterior of one simply connected compact set.

As in \cite{BLM}, let us show that this proposition implies the point (i) of Theorems~\ref{main2} and \ref{main1}.

\noindent {\bf Step 1: uniform $L^2$ estimate for $u^\eps-K_{\R^2}[\omega^\eps]$.}

For each time, we remark that $u^\eps - v^\eps[\omega^\eps]$ and $K_{\R^2}[\omega^\eps]-v^\eps[\omega^\eps]$ are divergence free, tend to zero when $|x|\to \infty$, have the same curl and circulations around $\K_{i,j}^\eps$ for all $i,j$. The only difference is that $u^\eps - v^\eps[\omega^\eps]$ is tangent to the boundary of $\Omega^\eps$, which implies that it is the Leray projection of $K_{\R^2}[\omega^\eps]-v^\eps[\omega^\eps]$. Therefore, by orthogonality of this projection in $L^2$ together with triangle inequality, we have 
\[
\| u^\eps-K_{\R^2}[\omega^\eps] \|_{L^2(\Omega^\eps)} \leq \| u^\eps-v^\eps[\omega^\eps] \|_{L^2(\Omega^\eps)} + \| v^\eps[\omega^\eps] -K_{\R^2}[\omega^\eps] \|_{L^2(\Omega^\eps)} \leq 2 \| v^\eps[\omega^\eps] -K_{\R^2}[\omega^\eps] \|_{L^2(\Omega^\eps)}.
\]
Under the assumptions of the point (i) of Theorems~\ref{main2} and \ref{main1}, the estimate of $\|\omega^\eps\|_{L^1\cap L^\infty}$ \eqref{om-cons} and Proposition~\ref{prop:key1} give that
\[
\| u^\eps-K_{\R^2}[\omega^\eps] \|_{L^2(\Omega^\eps)} \to 0 \text{ as }\eps\to 0 \quad \text{uniformly on time.}
\]
Recalling the standard estimate for the Biot-Savart kernel:
\begin{equation}\label{BS}
 \| K_{\R^2}[f] \|_{L^\infty(\R^2)}
\leq \Big\|\frac{1}{2\pi} \int_{\R^2}\frac{|f(y)|}{|x-y|} \, d y\Big\|_{L^\infty(\R^2)}
\leq C \|f\|_{L^1(\R^2)}^{1/2} \|f \|_{L^\infty(\R^2)}^{1/2},
\end{equation}
and the fact that $\| \text{\Large $\mathds{1}$}_{\R^2\setminus \Omega^\eps} \|_{L^2} \to 0$ (because $d_{\eps}/\eps \to \infty$ in the case of inclusions distributed on the square), we infer that 
\begin{equation}\label{conv:u-K}
 u^\eps-K_{\R^2}[\omega^\eps] \to 0 \quad \text{strongly in }L^\infty(\R^+;L^2(\R^2)), 
\end{equation}
where we have extended $u^\eps$ by zero inside the holes.

\noindent {\bf Step 2: compactness for the vorticity.}

Thanks to the uniform estimate of $\|\omega^\eps\|_{L^1\cap L^\infty}$ \eqref{om-cons}, Banach-Alaoglu's theorem infers that we can extract a subsequence such that
\[ \omega^\eps \rightharpoonup \omega\quad\text{ weak-$*$ in }L^\infty(\R^+; L^1\cap L^\infty(\R^2)),\]
which establishes the vorticity convergence stating in point (i) of Theorems~\ref{main2} and \ref{main1}.

Next, we derive a temporal estimate. For any $\varepsilon>0$, as $u^\eps$ is regular enough and tangent to the boundary, we can write the equation verified by $\omega^\eps$ for any test function $\phi \in H^1(\R^2)$:
\begin{eqnarray*}
(\partial_t \omega^\eps,\phi)_{H^{-1} , H^1} 
=\int_{\Omega^\eps} u^\eps \omega^\eps \cdot \nabla \phi = \int_{\R^2} (u^\eps- K_{\R^2}[\omega^\eps])\omega^\eps \cdot \nabla \phi + \int_{\R^2} K_{\R^2}[\omega^\eps]\omega^\eps \cdot \nabla \phi,
\end{eqnarray*}
which is bounded by $C\|\nabla \phi\|_{L^2}$ thanks to \eqref{om-cons}, \eqref{BS} and \eqref{conv:u-K}. Hence, we have
\[\| \partial_t \omega^\eps \|_{L^\infty(\R^+;H^{-1}(\R^2))} \leq C.\]
By Lemma C.1 in \cite{PLL}, this property can be used to extract a subsequence such that
\begin{equation}\label{conv:om}
 \omega^\eps \to \omega \text{ in } C([0,T]; L^{3/2}\cap L^4(\R^2)-w) \text{ for all }T.
\end{equation}

\noindent {\bf Step 3: compactness for the velocity.}

Now, we define $u:= K_{\R^2}[\omega]$ and we use the previous steps to pass to the limit in the decomposition
\begin{equation}\label{decomp2}
u^\eps-u = (u^\eps- K_{\R^2}[\omega^\eps]) + K_{\R^2}[\omega^\eps-\omega].
\end{equation}
Thanks to \eqref{conv:u-K}, it is obvious that the first term on the right-hand side of \eqref{decomp2} converges to zero in $L^2_{\loc}(\R^+ \times \R^2)$. Concerning the second term: for $x$ fixed, the map $y\mapsto \frac{(x-y)^\perp}{|x-y|^2}$ belongs to $L^{4/3}(B(x,1))\cap L^{3}(B(x,1)^c)$, then \eqref{conv:om} implies that for all $t,x$, we have 
\[
\int_{\R^2} \frac{(x-y)^\perp}{|x-y|^2} (\omega^\eps-\omega)(t,y)\, d y \to 0\quad \text{ as }\varepsilon\to 0.
\]
So, this integral converges pointwise to zero, and it is uniformly bounded by \eqref{BS} and \eqref{om-cons} with respect of $x$ and $t$. Applying the dominated convergence theorem, we obtain the convergence of $K_{\R^2}[\omega^\eps-\omega]$ in $L^2_{\loc}(\R^+ \times \R^2)$. This ends the proof of the velocity convergence stated in point (i) of Theorems~\ref{main2} and \ref{main1}. 

\noindent {\bf Step 4: passing to the limit in the Euler equations.}

Finally, we verify that $(u,\omega)$ is the unique solution of the Euler equations in $\R^2$.

The divergence and curl conditions are verified by the expression: $u=K_{\R^2}[\omega]$. Next, we use that $u^\eps$ and $\omega^\eps$ satisfies $\eqref{Euler5}$ in the sense of distribution in $\Omega^\eps$, that $u^\eps$ is regular and tangent to the boundary, to infer that for any test function $\phi\in C^\infty_c([0,\infty)\times \R^2)$, we have
\[
\int_0^\infty\int_{\R^2} \phi_t \omega^\eps \, d x d t+\int_0^\infty\int_{\R^2} \nabla\phi \cdot u^\eps \omega^\eps \, dx d t= -\int_{\R^2}\phi(0,x)\omega_0(x)\text{\Large $\mathds{1}$}_{\Omega^\eps} d x,
\]
where we have extended $\omega^\eps$ by zero and use that $\omega^\eps(0,\cdot)=\omega_{0}\vert_{\Omega^\eps}$.
By passing to the limit as $\varepsilon\to 0$, thanks to the strong-weak convergence of the pair $(u^\eps,\omega^\eps)$, we conclude that $(u,\omega)$ verifies the vorticity equation. In the whole plane, this is equivalent to state that $u$ verifies the velocity equation. As this solution is unique (Yudovich theorem), we deduce that the convergences hold without extracting a subsequence. This ends the proof of point (i) of Theorems~\ref{main2} and \ref{main1}.

\subsection{Convergence for small distance}\label{sec-cvg-small}

The proof of impermeability of the porous medium is totally different from the proof that the limit satisfies the non-linear equation in the vicinity of the inclusions.

 The easiest way to compute the limit of the tangency condition is to write the weak formulation. The weak form of 
\[
\div u_{0}= 0 \text{ in } \Omega \quad \text{and}\quad u_{0}\cdot n =0 \text{ on }\partial \Omega
\]
 reads
\begin{equation*}
\int_{\Omega} u_{0}\cdot h = 0 \text{ for all }h\in G(\Omega):=\{ w\in L^2(\Omega)\ : \ w=\nabla p, \text{ for some }p\in H^1_{\loc}(\Omega) \text{ and } w(x)=0 \text{ for large }x\}.
\end{equation*}
This general definition was introduced by Galdi \cite{Galdi} and used by G\'erard-Varet and Lacave \cite{GV-Lac2}. In this article we consider smooth test functions:
\[
h\in \{ \nabla p\ : \ p\in C^{\infty}_{c}(\overline{\Omega})\} 
\]
where
\begin{multline*}
C^{\infty}_{c}(\overline{\Omega}):=\{ \phi \text{ such that } D^k\phi \text{ is bounded and uniformly continuous in }\Omega, \ \forall k\in \N\\
\text{ and } \phi(x)=0 \text{ for large }x\}.
\end{multline*}

We refer to Galdi for discussions about the density of $\{ \nabla p\ : \ p\in C^{\infty}_{c}(\overline{\Omega})\}$ in $G(\Omega)$ for the $L^2$ norm (see e.g. in \cite{Galdi} the end of Section III.2 and Theorem II.7.2).

\begin{remark}
 Let us develop here an example, in order to focus on the importance of choosing test functions 
 which can have a jump across the segment. We consider $p\in C^{\infty}_{c}(\overline{\Omega_{1}})$ such that $p\equiv 0$ in a neighborhood of the endpoints and below the segment. Namely, there exists $\delta>0$ such that $p(x)=0$ for all $x\in B(0,\delta)\cup B((1,0),\delta)\cup ([0,1]\times[-\delta,0))$. Hence, assuming that $u_{0}$ is continuous up to the curve (with different values on each side), we have
\[
\int_{\Omega_{1}} u_{0} \cdot \nabla p= \int_{\delta}^{1-\delta} pu_{0}^{up}\cdot (-e_{2}).
\]
In this case, of course $\nabla p$ belongs to $\{ \nabla p\ : \ p\in C^{\infty}_{c}(\overline{\Omega})\}$, and we note that the tangency condition implies that the upper value of $u_{0}$ is tangent to the unit segment. It is important to note that if we only 
 consider smooth test functions in $\R^2$, we would only prove that $u_{0}^{up}\cdot e_{2} = u_{0}^{down}\cdot e_{2}$, which is weaker than the impermeability condition.
 
 We should not confuse the notation $C^{\infty}_{c}(\overline{\Omega})$ with the set of smooth functions compactly supported in the closure of $\Omega$. In particular, $C^{\infty}_{c}(\overline{\Omega_{1}})\subsetneq C^{\infty}_{c}(\R^2)$, since functions in $C^{\infty}_{c}(\overline{\Omega_{1}})$ are allowed to have a discontinuity though the segment.
\end{remark}
Therefore, we will prove the tangency condition in this sense:
\begin{equation}\label{tang weak}
\int_{\Omega} u_{0}\cdot \nabla p = 0 \text{ for all }p\in C^\infty_{c} (\overline{\Omega}).
\end{equation}
Similarly, the tangency condition on $u$ will read
\begin{equation}\label{tang u}
\int_{0}^\infty \int_{\Omega} u\cdot \nabla p = 0 \text{ for all }p\in C^{\infty}_{c}((0,\infty);C^\infty_{c} (\overline{\Omega})).
\end{equation}

The non-penetration of the porous medium for small distance will follow from the following key proposition.

\begin{proposition}\label{prop:key2} {\bf (Non-penetration)}
 Assume that $\K$ verifies {\rm (H1)} and let $(v^\eps)_{\eps}$ a bounded sequence in $L^2_{\loc}(\R^2)$, of divergence free vector fields which are tangent to the boundary $\partial\Omega^\eps$.
\begin{itemize}
 \item Consider $\Omega_{1}^\eps$ \eqref{Omega1eps} with $\K$ verifying also {\rm (H2)} and $\frac{d_{\eps}}{\eps^{2+\frac1\gamma}}\to 0$. For any $p\in C^\infty_{c} (\overline{\Omega_{1}})$, we have
\begin{equation*}
\int_{\Omega_{1}} v^\eps \cdot \nabla p \, dx \to 0\quad \text{as } \eps\to 0.
\end{equation*}
 \item Consider $\Omega_{2}^\eps$ \eqref{Omega2eps} (where $(\pm 1,0),(0,\pm 1)\in \partial\K$) and $\frac{d_{\eps}}\eps \to 0$ as $\eps\to 0$. For any $p\in C^\infty_{c} (\overline{\Omega_{2}})$, we have
\begin{equation*}
\int_{\Omega_{2}} v^\eps \cdot \nabla p \, dx \to 0\quad \text{as } \eps\to 0.
\end{equation*}
\end{itemize}
\end{proposition}

Unfortunately, we won't obtain an explicit rate in this proposition. The proof of this proposition will be performed in Section~\ref{sec:impermeability}.
In the rest of this section, we prove that this proposition implies the point (ii) of Theorems~\ref{main2} and \ref{main1}.

\noindent {\bf Step 1: weak convergence.}

By the estimates \eqref{est-utilde0}, \eqref{om-cons}, \eqref{estL2:ueps} and by Banach-Alaoglu's theorem, we can extract a subsequence $\varepsilon\to 0$ such that
\begin{eqnarray*}
&\omega^{\varepsilon} \rightharpoonup \omega& \text{ weak-$*$ in }L^\infty(\R^+;L^1\cap L^\infty(\Omega));\\
&u^\varepsilon_{0} \rightharpoonup u_0 &\text{ weakly in }L^2_{\loc}(\Omega);\\
&u^{\varepsilon} \rightharpoonup u &\text{ weak-$*$ in }L^\infty_{\loc}(\R^+;L^2_{\loc}(\Omega)).
\end{eqnarray*}
Moreover, it is obvious that:
\begin{equation*}
\omega^{\varepsilon}_{0}=\omega_{0}\vert_{\Omega^{\varepsilon}} \to \omega_{0} \text{ strongly in }L^q(\Omega),\quad \forall q\in[1,\infty).
\end{equation*}
In all these statements, we have extended the functions by zero inside the inclusions. Passing to the limit in the sense of distribution, it is straightforward that
\[
\div u=\div u_{0}=0 \text{ in }\Omega,\quad \curl u=\omega \text{ and } \curl u_{0}=\omega_{0} \text{ in }\Omega.
\]
The weak limit is also sufficient to prove that $\omega$ verifies \eqref{om-cons}.

\noindent {\bf Step 2: tangency condition.}

In both cases, we deduce directly from \eqref{est-utilde0} and Proposition~\ref{prop:key2} that for any $p\in C^{\infty}_{c}( \overline{\Omega})$ we have
\[
\int_{\Omega} u^\eps_{0}\cdot \nabla p \to 0.
\]
As $u^\eps_{0}$ converges weakly to $u_{0}$, we conclude that $u_{0}$ is tangent (in the weak sense \eqref{tang weak}).

Concerning $u$, for $p\in C^{\infty}_{c}((0,\infty)\times \overline{\Omega})$ we set $T,R>0$ such that $\supp p\subset(0,T) \times B(0,R)$. For a.e. $t\in [0,T]$, \eqref{estL2:ueps} and Proposition~\ref{prop:key2} imply that
\[
F(t):= \int_{\Omega} u^\eps(t,\cdot)\cdot \nabla p (t,\cdot) \to 0.
\]
Moreover, $F$ is uniformly bounded by $ \| u^\eps \|_{L^\infty(0,T;L^2( B(0,R)))} \| p \|_{L^\infty(\R^+;W^{1,\infty})}\leq C$ so the dominated convergence theorem together with the weak convergence of $u^\eps$ to $u$ give the tangency property for $u$ (in the weak sense \eqref{tang u}).

\begin{remark}\label{rem Grisvard}
Let us justify that the tangency is verified in a strong sense on the boundary, except on the corners. By uniqueness of vectors fields verifying the tangency in the weak sense, $\div u_{0}=0$ and $\curl u_{0} = \omega_{0}\in L^\infty$ (see \cite[Section 4.2]{GV-Lac2}), the standard elliptic estimates in domains with corners (the exterior of the segment can be seen as a domain with two corners of angles $2\pi$) imply that 
$u_{0}$ is continuous up to the boundary, except near the corners where it blows up (see e.g. \cite{Kondratiev,Mazya} and more recently \cite{Lac-uni}). Near the corners of the square, the velocity behaves like $1/| x |^{1/3}$. In the exterior of the segment, the velocity is also continuous up to the segment, with different values on each side of the segment, and blows up near the end-points like the inverse of the square root of the distance. 

This argument also holds for $u(t,\cdot)$ for a.e. $t$.
\end{remark}

\noindent {\bf Step 3: circulation condition.}

First, we justify that the circulation is well defined in the strong sense. As in Remark~\ref{rem Grisvard}, the standard elliptic estimates in domains with corners give the explicit behavior of $u_{0}$ and confirms the integrability of the trace. In the exterior of the segment, the quantity $\oint_{\partial \Omega_{1}} u_{0}\cdot \tau \, ds$ should be understood as the integral on the upper value $\int_{0}^1 u_{0}^{up} \cdot (-e_{1}) $ plus the lower value $\int_{0}^1 u_{0}^{down} \cdot e_{1}$. 

The easiest way to compute the limit of the circulation is to write the weak formulation. Namely, let us consider a smooth cutoff function $\chi$ such that $\chi \equiv 1$ in $B(0,R)$ and $\chi \equiv 0$ in $B(0,R+1)^c$, with $R$ large enough such that $\partial \Omega^\eps \subset B(0,R)$. Then, thanks to the regularity explained just above, we have:
\[
 \oint_{\partial \Omega} u_{0} \cdot \tau \, ds=-\int_{\Omega} \omega_{0} \chi\, dx - \int_{\Omega} u_{0} \cdot \nabla^\perp \chi\, dx.
\]
Moreover, we have by the Stokes formula and the zero circulation assumption
 \[
 -\int_{\Omega^\varepsilon} \omega_{0}^\varepsilon \chi\, dx - \int_{\Omega^\varepsilon} u_{0}^\varepsilon \cdot \nabla^\perp \chi\, dx = \sum_{i,j} \oint_{\partial \K_{i,j}^\varepsilon} u_{0}^\varepsilon \cdot \tau \, ds=0.
 \]
 Passing to the limit $\varepsilon\to 0$ in the previous equation gives that the circulation of $u_{0}$ around the unit segment is equal to zero.
 The proof for $u$ is exactly the same, adding a test function $\phi\in \mathcal{D}(\R^+)$ in time:
\[
\int_{\R+} \phi(t) \oint_{\partial \Omega} u(t,\cdot)\cdot \tau \, ds\, dt = 0.
\]

\noindent {\bf Step 4: passing to the limit in the non-linear equations.}

Of course, the weak convergences are not sufficient to pass to the limit in the non-linear term: $\int (u^\varepsilon \otimes u^\varepsilon):\nabla \phi$ for \eqref{Euler1} and $\int \omega^\varepsilon u^\varepsilon \cdot \nabla \psi$ for \eqref{Euler5}. We need some local strong convergence, far away from the porous medium. 
The following argument was introduced by Lions and Masmoudi \cite{LM99} and used 
 also in \cite{GV-Lac}.

Let us fix $T>0$ and $\mathcal{O}\Subset \Omega$ be a smooth open subset of $\Omega$. We write the Helmholtz-Weyl decomposition:
\[
u^\varepsilon = \mathbb{P}_{\mathcal{O}} u^\varepsilon + \nabla q^\varepsilon
\]
where $\mathbb{P}_{\mathcal{O}} u^\varepsilon$ is the Leray projection, i.e. a divergence free vector field which is tangent to the boundary $\partial \mathcal{O}$. Then, it is obvious that $\curl \mathbb{P}_{\mathcal{O}} u^\varepsilon= \curl u^\varepsilon=\omega^\varepsilon$ and $\Delta q^\varepsilon =0$. Moreover, by the orthogonality in $L^2$, we have:
\[
\| \mathbb{P}_{\mathcal{O}} u^\varepsilon \|_{L^\infty(0,T;L^2(\mathcal{O}))}+\| \nabla q^\varepsilon \|_{L^\infty(0,T;L^2(\mathcal{O}))} \leq 2 \| u^\varepsilon \|_{L^\infty(0,T;L^2(\mathcal{O}))} \leq C
\]
where we have used \eqref{estL2:ueps}. As $\mathbb{P}_{\mathcal{O}} u^\varepsilon $ is divergence free, bounded in $L^\infty(0,T;L^2(\mathcal{O}))$ and its curl belongs uniformly to $L^\infty(0,T\times \mathcal{O})$, it follows by standard elliptic estimates on smooth domains (e.g. Calder\'on-Zygmund inequality) that $\mathbb{P}_{\mathcal{O}} u^\varepsilon $ belongs uniformly to $L^\infty(0,T;H^1(\mathcal{O}))$. This argument can be also applied to state that $u^\varepsilon$ is uniformly bounded in $L^\infty(0,T;H^1(\mathcal{O}))$.

To derive a time estimate, let us denote
\[
\mathcal{V}(\mathcal{O}) = \text{ the closure in $H^{1}$ of } \{ \phi\in C^\infty_{c}(\mathcal{O}),\ \div \phi=0\},
\]
\[\mathcal{V}(\mathcal{O})' \text{ the dual of }\mathcal{V}(\mathcal{O})
\]
\begin{equation*}
\mathcal{H}(\mathcal{O}) = \text{ the closure in $L^2$ of } \{ \phi\in C^\infty_{c}(\mathcal{O}),\ \div \phi=0\}.
\end{equation*}
Hence we have that $\mathbb{P}_{\mathcal{O}} u^\varepsilon$ belongs uniformly to $L^\infty(0,T;\mathcal{H}(\mathcal{O}))$.
For any divergence free test function $\phi\in C^\infty_{c}((0,T)\times \mathcal{O})$, we compute thanks to \eqref{Euler1}:
\begin{eqnarray*}
\langle \partial_{t} \mathbb{P}_{\mathcal{O}} u^\varepsilon,\phi\rangle &=& \int_{0}^T \int_{\mathcal{O}} \mathbb{P}_{\mathcal{O}} u^\varepsilon \partial_{t} \phi\, dxdt = \int_{0}^T \int_{\mathcal{O}} u^\varepsilon \partial_{t} \phi\, dxdt = - \int_{0}^T \int_{\mathcal{O}} (u^\varepsilon\otimes u^\varepsilon):\nabla \phi\, dxdt\\
\Big| \langle \partial_{t} \mathbb{P}_{\mathcal{O}} u^\varepsilon,\phi\rangle \Big| &\leq& \| u^\varepsilon \|_{L^\infty(0,T;L^4(\mathcal{O}))}^2 \| \nabla \phi\|_{L^2(0,T;L^2(\mathcal{O}))}\sqrt{T} \leq C \| \phi\|_{L^2(0,T;\mathcal{V}(\mathcal{O}))}.
\end{eqnarray*}
In the last inequality, $C$ depends only on $\omega_{0}, \mathcal{O}, T$. By the Aubin-Lions lemma in $H^{1}(\mathcal{O}) \cap \mathcal{H}(\mathcal{O}) \hookrightarrow \mathcal{H}(\mathcal{O}) \hookrightarrow \mathcal{V}'(\mathcal{O})$ we conclude that there exists a subsequence such that
\[
\mathbb{P}_{\mathcal{O}} u^\varepsilon \to v \text{ strongly in } L^\infty(0,T;\mathcal{H}(\mathcal{O})).
\]
For this subsequence, we have that $\nabla q^\varepsilon=u^\varepsilon-\mathbb{P}_{\mathcal{O}} u^\varepsilon$ converges weakly in $L^2$ to a gradient $\nabla q$ (the set of gradient function is a closed subspace), hence at the limit we have $u=v+\nabla q$ and then $v=\mathbb{P}_{\mathcal{O}} u$. By the 
 unicity of the limit, we conclude that we do not need to extract a subsequence in the Aubin-Lions lemma:
\[
\mathbb{P}_{\mathcal{O}} u^\varepsilon \to \mathbb{P}_{\mathcal{O}} u \text{ strongly in } L^\infty(0,T;\mathcal{H}(\mathcal{O})),\text{ for the sequence $\varepsilon\to 0$ introduced in Step 1}.
\]

We are now in position to pass to the limit in the non linear equations.

For any divergence free test function $\phi\in C^{\infty}_{c}([0,\infty)\times \Omega)$, there exist $T>0$ and $\mathcal{O}\Subset \Omega$ be a smooth open subset of $\Omega$, such that $\supp \phi\subset[0,T)\times \mathcal{O}$. For $\varepsilon$ small enough such that $\mathcal{O}\Subset \Omega^\varepsilon$, the weak formulation of \eqref{Euler1} and \eqref{Euler4} reads:
\begin{equation*}
\int_{0}^\infty\int_{\Omega} u^\varepsilon \cdot \partial_{t} \phi\, dxdt+\int_{0}^\infty\int_{\Omega} (u^\varepsilon \otimes u^\varepsilon):\nabla \phi\, dxdt=-\int_{\Omega} u_{0}^\varepsilon\cdot \phi(0,\cdot) \, dx.
\end{equation*}
The weak limits of Step 1 are sufficient to pass to the limit in the first and third integrals. Before considering the second integral, let us note that for any harmonic function $f$ (i.e. $\Delta f=0$), we have the following relation:
\begin{equation}\label{algebra}
\int_{\mathcal{O}} \nabla f \otimes\nabla f : \nabla \phi =-\int_{\mathcal{O}} \div(\nabla f \otimes\nabla f) \phi= -\int_{\mathcal{O}}\left( \frac{1}{2} \nabla |\nabla f|^2 \cdot \phi + \Delta f \nabla f \cdot \phi \right)= 0,
\end{equation}
because $\phi$ is divergence free and compactly supported in $\mathcal{O}$. Such a relation can be applied with $f=q^\varepsilon$ and $f=q$. Therefore, the second integral of \eqref{weak Euler} can be decomposed as follows:
\begin{eqnarray*}
\int_{0}^\infty\int_{\Omega} (u^\varepsilon \otimes u^\varepsilon):\nabla \phi\, dxdt &=& \int_{0}^{T}\int_{\mathcal{O}} (\mathbb{P}_{\mathcal{O}}u^\varepsilon \otimes u^\varepsilon):\nabla \phi\, dxdt + \int_{0}^{T}\int_{\mathcal{O}} (\nabla q^\varepsilon \otimes \mathbb{P}_{\mathcal{O}} u^\varepsilon):\nabla \phi\, dxdt\\
&&+ \int_{0}^{T}\int_{\mathcal{O}} (\nabla q^\varepsilon \otimes \nabla q^\varepsilon ):\nabla \phi\, dxdt\\
&=&\int_{0}^{T}\int_{\mathcal{O}} (\mathbb{P}_{\mathcal{O}}u^\varepsilon \otimes u^\varepsilon):\nabla \phi\, dxdt + \int_{0}^{T}\int_{\mathcal{O}} (\nabla q^\varepsilon \otimes \mathbb{P}_{\mathcal{O}} u^\varepsilon):\nabla \phi\, dxdt
\end{eqnarray*}
thanks to \eqref{algebra}. By the weak limits of $u^\varepsilon$ and $\nabla q^\varepsilon$, the strong limit of $\mathbb{P}_{\mathcal{O}} u^\varepsilon$ and \eqref{algebra}, we easily verify that
\begin{align*}
\int_{0}^\infty\int_{\Omega} (u^\varepsilon \otimes u^\varepsilon):\nabla \phi\, dxdt \to&
\int_{0}^{T}\int_{\mathcal{O}} (\mathbb{P}_{\mathcal{O}}u \otimes u):\nabla \phi\, dxdt + \int_{0}^{T}\int_{\mathcal{O}} (\nabla q \otimes \mathbb{P}_{\mathcal{O}} u):\nabla \phi\, dxdt\\
&= \int_{0}^{T}\int_{\mathcal{O}} (\mathbb{P}_{\mathcal{O}}u \otimes u):\nabla \phi\, dxdt + \int_{0}^{T}\int_{\mathcal{O}} (\nabla q \otimes \mathbb{P}_{\mathcal{O}} u):\nabla \phi\, dxdt\\
&\ \ + \int_{0}^{T}\int_{\mathcal{O}} (\nabla q \otimes \nabla q ):\nabla \phi\, dxdt\\
&=\int_{0}^\infty\int_{\Omega} (u \otimes u):\nabla \phi\, dxdt.
\end{align*}
Therefore, we have proved that $u$ verifies the weak formulation of \eqref{Euler1} and \eqref{Euler4}:
\[
\int_{0}^\infty\int_{\Omega} u \cdot \partial_{t} \phi\, dxdt+\int_{0}^\infty\int_{\Omega} (u \otimes u):\nabla \phi\, dxdt=-\int_{\Omega} u_{0}\cdot \phi(0,\cdot) \, dx.
\]

The vorticity formulation is an obvious consequence of the velocity formulation: indeed, for any $\psi\in C^\infty_{c}([0,\infty)\times\Omega)$ we note that $\phi:=\nabla^\perp \psi$ is a divergence free test function for which the previous equation holds true. Thanks to elliptic regularity on a smooth subdomain $\mathcal{O}$ including the support of $\phi$, we have that $u$ belongs to $H^1(\mathcal{O})$ and then to $L^4(\mathcal{O})$. Hence, we can integrate by parts to get \eqref{Euler5}:
\[
\int_{0}^\infty\int_{\Omega} \omega \partial_{t} \psi\, dxdt+\int_{0}^\infty\int_{\Omega} \omega u \cdot \nabla \psi\, dxdt=-\int_{\Omega} \omega_{0} \psi(0,\cdot) \, dx.
\]
This ends the proof of point (ii) of Theorems~\ref{main2} and \ref{main1}.

\section{Permeability for large distance}\label{sec:permeability}
 
The goal of this section is to prove Proposition~\ref{prop:key1}. This section 
improves and simplifies the results of \cite{BLM}. 
The general argument is similar to \cite{BLM}, but we make more precise choices of 
cutoff functions. 
 We refer to that article for more details 
(in particular for the Biot-Savart law in exterior domains). Let $f\in L^\infty_{c}(\overline{\Omega^\eps})$ be fixed.

\subsection{Construction of the correction}

 The main idea is to use the explicit formula of the Green function (with Dirichlet boundary condition) in the exterior of one simply connected compact set $\K$:
\[
G_{\K}(x,y)=\frac1{2\pi} \ln \frac{| \Tc(x) - \Tc(y)|}{|\Tc(x)-\Tc(y)^*| |\Tc(y)|},
\]
where $\Tc: \ \K^c \to \R^2\setminus \overline{B(0,1)}$ is a biholomorphism such that 
\begin{equation}\label{T expansion}
\Tc(z)=\beta z +h(z)
\end{equation}
 for some $\beta \in \R^+$ and $h$ a bounded holomorphic function. Above, we have denoted by 
 $$y^*=\frac{y}{|y|^2}$$
 the conjugate point to $y$ across the unit circle in $\R^2$. Hence, it is verified in \cite[Section 3.1]{ILL} that the following vector field
\[
\nabla^\perp \int_{\R^2\setminus \K} G_{\K}(x,y) f(y) \, dy + \frac{\int_{\R^2\setminus \K} f}{2\pi} \nabla^\perp \ln | \Tc(x)|
\]
is divergence free, tangent to the boundary, goes to zero as $|x|\to \infty$, its curl is equal to $f$ and the circulation around $\K$ is equal to zero.

Now we introduce a cutoff function $\varphi_{i,j}^\eps$ equal to $1$ close to $\K_{i,j}^\eps$:
\begin{equation}\label{form:varphi}
 \varphi_{i,j}^\eps(x):= \varphi^\eps(x-z_{i,j}^\eps) 
\end{equation}
with $\varphi^\eps\in C^1$ such that $\varphi^\eps\equiv 1$ on $\frac\eps2 \partial\K$ and
\begin{itemize}
 \item $\varphi^\eps\equiv 0$ on $\partial ([-\tfrac{\eps+d_{\eps}}2,\tfrac{\eps+d_{\eps}}2]^2)$ if the inclusions are distributed in both directions;
 \item $\varphi^\eps\equiv 0$ on $\partial ([-\tfrac{\eps+d_{\eps}}2,\tfrac{\eps+d_{\eps}}2]\times [-\eps,\eps])$ if the inclusions are distributed only on the segment.
\end{itemize}

Then, the correction is defined by
\[
v^\eps[f]:= \nabla^\perp \psi^\eps,
\]
where 
\begin{align*}
\psi^\eps(x) :=& \frac1{2\pi} \Big(1- \sum_{i,j} \varphi_{i,j}^\eps(x) \Big) \int_{\Omega^\eps} \ln|x-y| f(y)\, d y\\
&+\frac{1}{2\pi}\sum_{i,j} \varphi_{i,j}^\eps (x)\int_{\Omega^\eps}{\ln}\frac{ \eps |\Tca{i,j}(x)-\Tca{i,j}(y)||\Tca{i,j}(x)|}{2\beta|\Tca{i,j}(x)-\Tca{i,j}(y)^*|} f(y) \, d y, 
\end{align*}
with
\[
\Tca{i,j}(x) := \Tc\left(\frac{x-z_{i,j}^\eps}{\eps/2} \right)\ : \ (\K^\eps_{i,j})^c \to \R^2\setminus \overline{B(0,1)}.
\]
In the neighborhood of $\K_{i,j}^\eps$, this correction corresponds to the Biot-Savart law in the exterior of one obstacle, whereas, far away the porous medium, it is equal to the Biot Savart law in the whole plane $\R^2$. More precisely, we can check that $v^\eps[f]$ verifies 
the following properties:
\begin{equation*}
\div v^\eps[f] =0 \text{ in } \Omega^\varepsilon, \quad v^\eps[f] \cdot n =0 \text{ on } \partial \Omega^\varepsilon, \quad \lim_{x\to \infty} |v^\eps[f](x)|=0.
\end{equation*}

We decompose $K_{\R^2}[f]-v^\eps[f]$ as 
\begin{equation} \label{decompo we}
K_{\R^2}[f]-v^\eps[f] =\frac1{2\pi} \sum_{i,j} \nabla^\perp \varphi_{i,j}^\eps(x) (w_{i,j}^{1,\eps}+w_{i,j}^{2,\eps}) + \varphi_{i,j}^\eps(x) (w_{i,j}^{3,\eps}+w_{i,j}^{4,\eps}) ,
\end{equation}
where
\begin{equation*}\begin{split}
w_{i,j}^{1,\eps}(x)=& \int_{\Omega^\eps} \ln \frac{2\beta|x-y|}{\eps|\Tca{i,j}(x)-\Tca{i,j}(y)|}f(y)\, dy, \\
w_{i,j}^{2,\eps}(x)=& \int_{\Omega^\eps} \ln \frac{|\Tca{i,j}(x)-\Tca{i,j}(y)^*|}{|\Tca{i,j}(x)|}f(y)\, dy, \\
w_{i,j}^{3,\eps}(x)=& \int_{\Omega^\eps} \Biggl(\frac{(x-y)^\perp}{|x-y|^2}- (D\Tca{i,j})^T(x)\frac{(\Tca{i,j}(x)-\Tca{i,j}(y))^\perp}{|\Tca{i,j}(x)-\Tca{i,j}(y)|^2} \Biggl) f(y)\, dy, \\
w_{i,j}^{4,\eps}(x)=& (D\Tca{i,j})^T(x) \int_{\Omega^\eps} \Biggl(\frac{\Tca{i,j}(x)-\Tca{i,j}(y)^*}{|\Tca{i,j}(x)-\Tca{i,j}(y)^*|^2}- \frac{\Tca{i,j}(x)}{|\Tca{i,j}(x)|^2}\Biggl)^\perp f(y)\, dy.
\end{split}\end{equation*}

In the following subsection, we estimate $w_{i,j}^{k,\eps}$ on the support of $\varphi_{i,j}^\eps$, and next, we will look for the best cutoff function $\varphi^\eps$.

\subsection{Cell problem estimates}

When $\K=\overline{B(0,1)}$, $\Tc={\rm Id}$ (so $\beta=1$) and $w_{1}^\eps=w_{3}^\eps \equiv 0$. In this case, we also have $\Tca{i,j}(x)-\Tca{i,j}(y)^* =\dfrac2\eps\Big(x-z_{i,j}^\eps -\eps^2 \dfrac{y-z_{i,j}^\eps}{4|y-z_{i,j}^\eps|^2}\Big)$. Except in an $\varepsilon$-neighborhood of the inclusion, we note that $\Tca{i,j}(y)^*$ is small compare to $\Tca{i,j}(x)$. Hence, we can guess that $w_{2}^\eps$ and $w_{4}^\eps$ are small. This remark is the main motivation of this decomposition, and in the following estimates, we split the integrals in two parts: a small area in the vicinity of the inclusion and the far away region where $\Tc$ behaves as $\beta\ {\rm Id}$. These estimates were made in \cite{BLM}, but we include here the main arguments for the sake of completeness. Moreover, we slightly improve the estimate of $w_{i,j}^{1,\eps}$. 

Let us before recall some basic estimates on conformal mapping.
From the definition of $ \Tca{i,j}$, it is possible to get the following Lipschitz estimates (see Lemma 2.1 in \cite{BLM}):
\begin{equation}\label{Lip}
\| \Tca{i,j} \|_{\mathrm{Lip}} \leq \frac{C}{\varepsilon} \quad \text{ and } \quad \| (\Tca{i,j})^{-1} \|_{\mathrm{Lip}} \leq C \varepsilon.
\end{equation}
Moreover, as $\Tc$ behaves at infinity as $\beta\ {\rm Id}$, it is also natural (Lemma 2.2 in \cite{BLM}) that for all $r>0$
\begin{equation}\label{anneau1}
 \Tca{i,j}\Bigl(\partial B(z^\eps_{i,j},r)\cap (\K^\eps_{i,j})^c\Bigr) \subset B\Bigl(0,C_{1}\frac r\eps\Bigr)\setminus B\Bigl(0,C_{2}\frac r\eps\Bigr)
\end{equation}
and
\begin{equation}\label{anneau2}
 (\Tca{i,j})^{-1}\Bigl(\partial B(0,r+1)\Bigr) \subset B\Bigl(z^\eps_{i,j},\eps C_{3}(r+1)\Bigr)\setminus B\Bigl(z^\eps_{i,j},\eps C_{4}(r+1)\Bigr).
\end{equation}
for some $C_{1},C_{2},C_{3},C_{4}$ positive numbers independent of $i,j, \eps$.

\noindent {\bf Estimate of $w_{i,j}^{1,\eps}$.}

For $x\in\supp \varphi^\eps_{i,j}$ fixed, we decompose the integral in two parts: 
\begin{equation}\label{eq.loinpres}
\begin{split}\Omega^\eps_{C}:=\{y\in \Omega^\eps,\ |\Tca{i,j}(x)-\Tca{i,j}(y)|\leq \eps^{-1/4} \},\\
\Omega^\eps_{F}:=\{y\in \Omega^\eps,\ |\Tca{i,j}(x)-\Tca{i,j}(y)|> \eps^{-1/4} \}.
\end{split}
\end{equation}
In the subdomain close to the inclusion $\Omega^\eps_{C}$, we set $z=\eps\Tca{i,j}(x)$ and we change variables $\eta=\eps \Tca{i,j}(y)$:
\begin{equation*}\begin{split}
\int_{\Omega^\eps_{C}} \Bigl| \ln (\eps|\Tca{i,j}(x)-\Tca{i,j}(y)|) f(y)\Bigl| \, d y
&\leq \int_{B(z, \eps^{3/4})} \Bigl| \ln |z-\eta| f((\Tca{i,j})^{-1}(\tfrac\eta\eps))\Bigl| \frac{ \bigl| \det D(\Tca{i,j})^{-1}\bigr|(\tfrac\eta\eps)}{\varepsilon^2}\, d \eta\\
&\leq \int_{B(z,\eps^{3/4})} \Bigl| \ln |z-\eta| f((\Tca{i,j})^{-1}(\tfrac\eta\eps))\Bigl|\tfrac14 \bigl| \det D\Tc^{-1}\bigr|(\tfrac\eta\eps) d \eta.
\end{split}\end{equation*}
Using that $D\Tc^{-1}$ and $f$ are bounded functions, we compute that: 
\begin{equation*}
\int_{\Omega^\eps_{C}} \Bigl| \ln (\eps|\Tca{i,j}(x)-\Tca{i,j}(y)|) f(y)\Bigl| \, d y
\leq C\| f\|_{L^\infty} \int_{B(0, \eps^{3/4} )} \Bigl| \ln |\xi| \Bigl| \, d \xi \leq C \| f\|_{L^\infty} \eps^{3/2}|\ln \eps| .
\end{equation*}
To deal with $\ln(2\beta|x-y|)$, we remark that if $y\in \Omega^\eps_{C}$, then \eqref{Lip} gives
\[
|x-y|=|(\Tca{i,j})^{-1}(\Tca{i,j}(x)) - (\Tca{i,j})^{-1}(\Tca{i,j}(y))| \leq \eps C |\Tca{i,j}(x)- \Tca{i,j}(y)| \leq C \eps^{3/4}.
\]
So, we have
\begin{align*}
\int_{\Omega^\eps_{C}} \Bigl| \ln( {2\beta|x-y|})f(y) \Bigl|\, d y 
&\leq \int_{B(x, C\eps^{3/4})} \Bigl| \ln (2\beta|x-y|) f(y) \Bigl|\, d y \\
&\leq \|f\|_{L^\infty} \int_{B(0, C\eps^{3/4})} \Bigl| \ln |2\beta \xi| \Bigl|\, d\xi 
\leq C\|f\|_{L^\infty} \eps^{3/2} |\ln \eps|.
\end{align*}

In the subdomain far away from the inclusion $\Omega^\eps_{F}$, we have by \eqref{Lip}
\[
\eps^{-1/4} \leq |\Tca{i,j}(x) - \Tca{i,j}(y)| \leq \eps^{-1} C |x-y|,
\]
hence $|x-y|\geq \frac{ \eps^{3/4}}C$. 
Therefore, with $h$ defined in \eqref{T expansion}, writing
\begin{equation}\label{eq.ln2}
\ln \frac{\eps|\Tca{i,j}(x)-\Tca{i,j}(y)|}{2\beta|x-y|}
= \ln \frac{\Bigl|2\beta(x-y) + \eps \Bigl(h\big(\frac{x-z^\eps_{i,j}}{\eps/2}\big)- h\big(\frac{y-z^\eps_{i,j}}{\eps/2}\big)\Bigr)\Bigr|}{2\beta|x-y|},
\end{equation}
we have
\[
\frac{\eps \Bigl| h\big(\frac{x-z^\eps_{i,j}}{\eps/2}\big)- h\big(\frac{y-z^\eps_{i,j}}{\eps/2}\big)\Bigl| }{2\beta|x-y|} 
\leq \frac{\|h\|_{L^\infty}}{\beta } C \eps^{1/4},
\]
which is smaller that $1/2$ for $\varepsilon$ small enough.
We note easily that
\begin{equation}\label{est ln}
\Bigl|\ln\tfrac{|b+c|}{|b|} \Bigl| \leq 2 \tfrac{|c|}{|b|},\qquad\mbox{ if }\tfrac{|c|}{|b|}\leq \tfrac12.
\end{equation}
Applying this inequality with $c=h\big(\frac{x-z^\eps_{i,j}}{\eps/2}\big)- h\big(\frac{y-z^\eps_{i,j}}{\eps/2}\big)$ and $b= \tfrac{2\beta(x-y)}\eps$, we compute from \eqref{eq.ln2}:
\begin{eqnarray*}
\Biggl|\ln \frac{\eps|\Tca{i,j}(x)-\Tca{i,j}(y)|}{\beta|x-y|}\Biggl|
&\leq& 2\frac{\eps|h(\frac{x-z^\eps_{i,j}}{\eps/2})-h(\frac{y-z^\eps_{i,j}}{\eps/2})|}{2\beta|x-y|} \leq \frac{2\eps \|h\|_{L^\infty}}{\beta|x-y|}.
\end{eqnarray*}
Therefore, using \eqref{BS}, we obtain
\begin{eqnarray*}
\int_{\Omega^\eps_{F}} \Bigl|\ln \frac{\beta|x-y|}{\eps|\Tca{i,j}(x)-\Tca{i,j}(y)|} f(y) \Bigl| \, d y
&\leq &\int_{ \Omega^\eps}\frac{2\eps \|h\|_{L^\infty}}{\beta|x-y|} |f(y)|\, d y\\
&\leq &C\eps \| f\|_{L^\infty}^{1/2} \| f\|_{L^1}^{1/2} 
\end{eqnarray*}
which allows us to conclude that
\begin{equation}\label{est:w1ij}
 \| w_{i,j}^{1,\eps} \|_{L^\infty (\supp \varphi^\eps_{i,j})} \leq C\varepsilon \| f\|_{L^1\cap L^\infty} 
\end{equation}
with $C$ independent of $i,j, \eps$ and $f$.

\noindent {\bf Estimate of $w_{i,j}^{2,\eps}$.}

We set $z=\eps \Tca{i,j}(x)$, and changing variables $\eta = \eps\Tca{i,j}(y)$, we get
\begin{equation*}
w_{i,j}^{2,\eps}(x)= \int_{B(0,\eps)^c} \ln \frac{|z- \eps^2 \eta^*|}{|z|}f(\tfrac\eps2 \Tc^{-1}(\tfrac\eta\eps)+z^\eps_{i,j}) \tfrac14 |\det D\Tc^{-1}|(\tfrac\eta\eps) \, d \eta.
\end{equation*}
Without loss of generality, we assume that $0\in \stackrel{\circ}{\K}$, and we fix $\delta >0$ so that $B(0,\delta)\subset \K$. For $x\in \supp \varphi^\eps_{i,j}$, we easily note that 
\[
\delta \eps \leq |x-z^\eps_{i,j}| \leq \sqrt{2}\frac{\eps+d_{\eps}}2,
\]
then, we deduce by \eqref{anneau1} that
\[
C_{2} \delta \eps \leq |z| \leq C_{1}\sqrt{2}\frac{\eps+d_{\eps}}2.
\]
As for any $\eta$ we have
\begin{equation*}
\frac{|\eps^2 \eta^*|}{|z|} \leq \frac{\eps}{C_{2}\delta |\eta|},
\end{equation*}
we infer by \eqref{est ln} (with $b=z$ and $c=-\eps^2 \eta^*$) that
\begin{equation*}
\left| \ln \frac{|z- \eps^2 \eta^*|}{|z|} \right|
\leq 2 \frac{\eps^2 |\eta^*|}{|z|}
\leq \frac{2\eps}{C_{2} \delta |\eta|}\qquad \mbox{ if }\quad \frac{\eps}{C_{2}\delta |\eta|}\leq \frac12.
\end{equation*}

Therefore, we define $R=2/(C_{2}\delta)$ and we split the integral in two parts: $ B(0,R\eps)^c$ and $B(0,R\eps)\setminus B(0,\eps)$. 

In the first subdomain $B(0,R\eps)^c$, we use the previous inequality to compute
\begin{equation*}\begin{split}
\Bigl|\int_{B(0,R\eps)^c}& \ln \frac{|z- \eps^2 \eta^*|}{|z|}f(\tfrac\eps2 \Tc^{-1}(\tfrac\eta\eps)+z^\eps_{i,j})\tfrac14 |\det D\Tc^{-1}|(\tfrac\eta\eps) \, d \eta\Bigl|\\
&\leq \frac{2\eps}{C_{2}\delta} \int_{\R^2} \frac{|f(\tfrac\eps2 \Tc^{-1}(\tfrac\eta\eps)+z^\eps_{i,j})| \tfrac14 |\det D\Tc^{-1}|(\tfrac\eta\eps)}{|\eta|} \, d \eta\\
&\leq C\eps \Big\| f(\tfrac\eps2 \Tc^{-1}(\tfrac\eta\eps)+z^\eps_{i,j}) \tfrac14 \det D\Tc^{-1}(\tfrac\eta\eps) \Big\|_{L^\infty}^{1/2}
\Big\| f(\tfrac\eps2 \Tc^{-1}(\tfrac\eta\eps)+z^\eps_{i,j}) \tfrac14\det D\Tc^{-1}(\tfrac\eta\eps) \Big\|_{L^1}^{1/2}\\
&\leq C\eps \|f\|_{L^\infty}^{1/2}\|f\|_{L^1}^{1/2},
\end{split}\end{equation*}
where we have applied \eqref{BS} for the function $\eta \mapsto |f(\tfrac\eps2 \Tc^{-1}(\tfrac\eta\eps)+z^\eps_{i,j})|\tfrac14 |\det D\Tc^{-1}|(\tfrac\eta\eps)$ at $x=0$, used that $D \Tc^{-1}$ is bounded and that $\| f(\tfrac\eps2 \Tc^{-1}(\tfrac\eta\eps)+z^\eps_{i,j})\tfrac14 \det D\Tc^{-1}(\tfrac\eta\eps) \|_{L^1}=\|f\|_{L^1}$ by changing variables back.

In the second subdomain $B(0,R\eps)\setminus B(0,\eps)$, we use the relation
\[
\frac{|z- \eps^2 \eta^*|}{|z|}=\frac{|\eta- \eps^2 z^*|}{|\eta|}
\]
which can be easily verified by squaring both side. As $D\Tc^{-1}$ is bounded and $\eps^2 z^*\in B(0,\eps)$, we compute
\begin{equation*}\begin{split}
\Bigl|\int_{B(0,R\eps)\setminus B(0,\eps)} \ln \frac{|z- \eps^2 \eta^*|}{|z|}f(\tfrac\eps2 \Tc^{-1}(\tfrac\eta\eps)+z^\eps_{i,j}) \tfrac14 |\det D\Tc^{-1}|(\tfrac\eta\eps) \, d \eta\Bigl|
&\leq 2C \| f\|_{L^\infty}\int_{B(0,(R+1)\eps)} | \ln |\eta| |\, d\eta\\
&\leq C\| f\|_{L^\infty} \eps^2 |\ln \eps|,
\end{split}\end{equation*}
which allows us to conclude that
\begin{equation}\label{est:w2ij}
 \| w_{i,j}^{2,\eps} \|_{L^\infty (\supp \varphi^\eps_{i,j})} \leq C\varepsilon \| f\|_{L^1\cap L^\infty} 
\end{equation}
with $C$ independent of $i,j, \eps$ and $f$.

\noindent {\bf Estimate of $w_{i,j}^{3,\eps}$ and $w_{i,j}^{4,\eps}$.}

With similar technics, i.e. by changing variables and using the expression of $\Tca{i,j}$ in terms of $\Tc$, it is not difficult to prove that 
\begin{equation}\label{est:w34ij}
 \| w_{i,j}^{3,\eps} \|_{L^\infty (\Omega^\eps)}+ \| w_{i,j}^{4,\eps} \|_{L^\infty (\Omega^\eps)} \leq C \| f\|_{L^1\cap L^\infty}, 
\end{equation}
with $C$ independent of $i,j, \eps$ and $f$ (see e.g. \cite[Theorem 4.1]{ILL}). Actually, we were more precise in \cite{BLM} proving an estimate which tends to zero, but as we will discuss in the next subsection, the restriction $\frac{d_{\eps}}{\eps}\to \infty$ (for inclusions distributed on the square) comes from $w_{i,j}^{1,\eps}$ and $w_{i,j}^{2,\eps}$.

\subsection{Optimal cutoff function}\label{subsec:cutoff}

Putting together the form of $\varphi^\eps_{i,j}$ \eqref{form:varphi}, the decomposition \eqref{decompo we} and the estimates \eqref{est:w1ij}, \eqref{est:w2ij} and \eqref{est:w34ij} we have obtained:
\begin{equation}\label{est:permeability}
 \| K_{\R^2}[f]-v^\eps \|_{L^2(\Omega^\eps)} \leq C \|f\|_{L^1\cap L^\infty} \Bigl(\eps \| \nabla \varphi^\eps \|_{L^2} + \| \varphi^\eps \|_{L^2} \Bigl) \sqrt{\sum_{i,j} 1},
\end{equation}
so the question is to find the best $\varphi^\eps$ such that the right hand side term tends to zero. 

\noindent {\bf Case 1: inclusions distributed on the square.}

In the case \eqref{Omega2eps}, the number of inclusions is $N_{\eps}^2 \leq 2/d_{\eps}^2$, so $\sqrt{\sum_{i,j} 1}\leq \sqrt 2/d_{\eps}$.

For $d_{\eps} \geq \eps$, it is sufficient to consider $\varphi^\eps$ independent of $d_{\eps}$. Namely, let $\varphi^1 \in C^\infty_{c}([-1,1]^2)$ such that $\varphi^1\equiv 1$ on $[-1/2,1/2]^2$, then we define
\[
\varphi^\eps (x):= \varphi^1\Big(\frac x\eps\Big).
\]
It is easy to note that $\varphi^\eps \equiv 1$ on $\frac\eps2\partial\K$ and $\varphi^\eps\equiv 0$ on $\partial ([-\tfrac{\eps+d_{\eps}}2,\tfrac{\eps+d_{\eps}}2]^2)$. With this expression, it is also clear that there exists $C>0$ such that
\[
 \| \varphi^\eps \|_{L^2} ,\ \eps \| \nabla \varphi^\eps \|_{L^2} \leq C\eps.
\]

In this case, \eqref{est:permeability} reads
\begin{align*}
 \| K_{\R^2}[f]-v^\eps \|_{L^2(\Omega^\eps)}\leq C \|f\|_{L^1\cap L^\infty} \frac\eps {d_{\eps}}
\end{align*}
which ends the proof of the first case of Proposition~\ref{prop:key1}.

\noindent {\bf Case 2: Inclusions distributed on the segment.}

We now have 
\begin{equation}\label{miniOFN}
 \sqrt{\sum_{i,j} 1} =\sqrt{N_{\eps}} \leq \min\Big(\sqrt{\frac2\eps},\sqrt{\frac2{d_\eps}}\Big). 
\end{equation}
In \cite{BLM}, we have considered a cutoff function such that $\supp \varphi^\eps$ was included in $([-\tfrac{\eps+d_{\eps}}2,\tfrac{\eps+d_{\eps}}2]^2) \setminus ([-\tfrac{\eps}2,\tfrac{\eps}2]^2)$, which allowed us to get a zero limit for $d_{\eps}/\eps^2\to \infty$. This rate is optimal if $\K$ is a square, but for the other shape, we can improve the cutoff function by using $\gamma$ (see (H2)).

The assumption (H2) implies that the minimal distance between $\K^\varepsilon_{i,1}$ and $\K^\varepsilon_{i+1,1}$ is reached for a unique point of $\partial \K^\varepsilon_{i,1}$: $z^\varepsilon_{i,1}+(\frac\eps2,0)$. Up to choosing $\rho_{2}$ smaller, we can assume without loss of generality that 
\begin{equation*}
1-2^{\gamma+1} \rho_{2}>0.
\end{equation*}

The idea is to define a cutoff function which depends on the space between $\frac\eps2\K^\varepsilon$ and the line $\{(\tfrac{\eps+d_{\eps}}2, \R)\}$. Hence, we define
\[
d(x_{2})=\tfrac{d_{\eps}}2+\tfrac{\varepsilon}2 \rho_{2} |2x_{2}/\varepsilon|^{\gamma+1}
\]
and we remark from (H2) that
\[
\Big( [-\tfrac{\eps+d_{\eps}}2 , -\tfrac{\eps+d_{\eps}}2 + d(x_{2})] \cup [\tfrac{\eps+d_{\eps}}2-d(x_{2}) , \tfrac{\eps+d_{\eps}}2 ] \Big) \times\{x_{2}\} \subset \Big([-\tfrac{\eps+d_{\eps}}2, \tfrac{\eps+d_{\eps}}2] \times [-\eps,\eps] \Big)\setminus \tfrac\eps2 \K \text{ for all } x_{2}\in [-\varepsilon,\varepsilon] .
\]
Let $\varphi\in \mathcal{C}^\infty(\R)$ be a positive non-increasing function such that $\varphi(s)=1$ if $s\leq 0$ and $\varphi(s)=0$ if $s\geq 1$. Then we introduce:
\[
\varphi^\varepsilon(x)=\varphi\Big(\frac{2|x_{2}|-\varepsilon}\varepsilon\Big)\Big[ 1-\varphi\Big(\frac{(\varepsilon+d_{\eps})/2-x_{1}}{d(x_{2})}\Big)-\varphi\Big(\frac{x_{1}+(\eps+d_{\eps})/2}{d(x_{2})}\Big)\Big]
\]
We can check that $\varphi^\varepsilon(x)\in [0,1]$ for all $x\in \R^2$, that $\varphi^\varepsilon(x)=0$ for all $x\in ([-\tfrac{\eps+d_{\eps}}2,\tfrac{\eps+d_{\eps}}2]\times [-\varepsilon,\varepsilon])^c$ and that $\varphi^\varepsilon(x)=1$ in the neighborhood of $\frac\eps2 \K$. Hence,
\[
\| \varphi^\varepsilon \|_{L^2}\leq \sqrt{2 \eps(\eps + d_{\eps})}.
\]

Next, we see that for all $x$ we have
\[
|\nabla\varphi^\varepsilon(x)| \leq \frac{C}\varepsilon + \frac{C}{d(x_{2})}+ \frac{Cd(x_{2})}{d(x_{2})^2}\leq \frac{C}\varepsilon + \frac{C}{d(x_{2})},
\]
where we have used that on the support of $\varphi'\Big(\frac{(\varepsilon+d_{\eps})/2-x_{1}}{d(x_{2})}\Big)$ we have $|(\varepsilon+d_{\eps})/2-x_{1}|\leq d(x_{2})$. So we compute:
\begin{equation*}
\begin{split}
\| \nabla \varphi^\varepsilon \|_{L^2}\leq& C\Big(\frac{2\eps(\varepsilon+d_{\eps})}{\varepsilon^2}+ 4\int_{0}^{\varepsilon}\int_{0}^{d(x_{2})} \frac1{d(x_{2})^2} \, dx_{1}dx_{2} \Big)^{1/2}\leq C\Big(1+\frac{d_{\eps}}{\eps}+ 4\int_{0}^{\varepsilon} \frac1{d(x_{2})} \, dx_{2} \Big)^{1/2}\\
& \leq C\Big(1+\frac{d_{\eps}}{\eps}+ 8\int_{0}^{\varepsilon} \frac1{d_{\eps}+\varepsilon^{-\gamma} \rho_{2} (2x_{2})^{\gamma+1}} \, dx_{2} \Big)^{1/2}\\
&\leq C\Big(1+\frac{d_{\eps}}{\eps}+ 8\int_{0}^\delta \frac1{d_{\eps}} \, dx_{2} + 2^{2-\gamma}\int_{\delta}^{\varepsilon} \frac1{\varepsilon^{-\gamma} \rho_{2} x_{2}^{\gamma+1}} \, dx_{2} \Big)^{1/2}\\
&\leq C\Big(1+\frac{d_{\eps}}{\eps}+ 8\frac{\delta}{d_{\eps}} + \frac{2^{2-\gamma}}{\gamma\varepsilon^{-\gamma} \rho_{2} \delta^{\gamma}} \Big)^{1/2}
\end{split}
\end{equation*}
where we choose $\delta \in (0,\eps)$ such that $\frac{\delta}{d_{\eps}}\leq \frac1{\varepsilon^{-\gamma} \delta^{\gamma}}$, i.e. $\delta=\min(\eps, (d_{\eps} \eps^\gamma)^{1/(\gamma+1)})$. 

In the case where $d_{\eps} \geq \eps$, then we take $\delta =\eps$ in the previous computation and we get 
\[
\| \nabla \varphi^\varepsilon \|_{L^2}\leq C\Big(9 +\frac{d_{\eps}}{\eps} \Big)^{1/2}\leq C + C\Big(\frac{d_{\eps}}{\eps}\Big)^{1/2}
\]
which implies by \eqref{est:permeability} and \eqref{miniOFN} that
\[
\| K_{\R^2}[f]-v^\eps \|_{L^2(\Omega^\eps)} \leq C \|f\|_{L^1\cap L^\infty} \Bigl(\frac{C\eps}{\sqrt{\eps}}+\frac{C\eps}{\sqrt{d_\eps}}\Big(\frac{d_{\eps}}{\eps}\Big)^{1/2} + \frac{\sqrt{2 \eps(\eps + d_{\eps})}}{\sqrt{d_{\eps}}}\Big)\leq C \|f\|_{L^1\cap L^\infty} \sqrt{\eps}.
\]

In the other case $d_{\eps} < \eps$, then we take $\delta= (d_{\eps} \eps^\gamma)^{1/(\gamma+1)}$ which gives
\[
\| \nabla \varphi^\varepsilon \|_{L^2}\leq C \Big(1 + \frac{(d_{\eps} \eps^\gamma)^{1/(\gamma+1)}}{d_{\eps}}\Big)^{1/2}\leq C \Big(1 + \Big(\frac{\eps}{d_{\eps}}\Big)^{\frac{\gamma}{2(\gamma+1)}}\Big).
\]
Together with \eqref{est:permeability} and \eqref{miniOFN}, we conclude that in this case
\begin{align*}
 \| K_{\R^2}[f]-v^\eps \|_{L^2(\Omega^\eps)} &\leq C \|f\|_{L^1\cap L^\infty} \Bigl(\frac{C\eps}{\sqrt{\eps}}+\frac{C\eps}{\sqrt{\eps}}\Big(\frac{\eps}{d_{\eps}}\Big)^{\frac{\gamma}{2(\gamma+1)}} + \frac{\sqrt{2 \eps(\eps + d_{\eps})}}{\sqrt{\eps}}\Big)\\
 &\leq C \|f\|_{L^1\cap L^\infty} \sqrt{\eps} \Bigl(1+ \Big(\frac{\eps}{d_{\eps}}\Big)^{\frac{\gamma}{2(\gamma+1)}} \Big).
\end{align*}

This ends the proof of Proposition~\ref{prop:key1}.

\begin{remark}\label{rem:cutoff}
For the Laplace, Stokes and Navier-Stokes problems, we do not have $\eps$ in front of $ \| \sum_{i,j}\nabla \varphi^\eps \|_{L^2}$. Finding $\varphi^\eps$ which minimizes $\| \nabla \varphi^\eps\|_{L^2}$ is a standard problem in optimization, and it is well known that the minimizer is the solution of $\Delta \varphi^\eps =0$. In an annulus $B(0,d_{\eps})\setminus B(0,\eps)$, the solution is 
\[
\frac{\ln \frac{|x|}{d_{\eps}}}{\ln \frac \eps {d_{\eps}}}
\]
and we get in the case of inclusions distributed in the square, for $d_{\eps}\gg \eps$,
\[
\| \sum_{i,j}\nabla \varphi^\eps \|_{L^2} \leq \frac{1}{d_{\eps}}\Big(\ln \frac {d_{\eps}}{\eps } \Big)^{-1/2} \leq \frac C{d_{\eps} |\ln \eps|^{1/2}},
\]
so the behavior depends on the limit of $d_{\eps} |\ln \eps|^{1/2}$ (see \cite{Allaire90a,Allaire90b,CM82}).

When the inclusions are distributed on the segment, the number of obstacles is $N_{\eps}\leq C/d_{\eps}$ and then
\[
\| \sum_{i,j}\nabla \varphi^\eps \|_{L^2} \leq \frac{1}{\sqrt{d_{\eps}}}\Big(\ln \frac {d_{\eps}}{\eps } \Big)^{-1/2} \leq \Big(\frac C{d_{\eps} |\ln \eps|} \Big)^{1/2} ,
\]
where the limit $d_{\eps} |\ln \eps|$is the criterion in \cite{Allaire90b}.

For inviscid flow, thanks to the additional $\eps$ in \eqref{est:permeability}, we manage to treat smaller distance than for viscous fluid. When the inclusions are distributed on the segment, we even deal with distance smaller than the inclusion size.
\end{remark}

\begin{remark}
To obtain a correction, the basic idea could be to cutoff the stream function and then to define:
\[
\tilde v^\varepsilon[f](x):= \nabla^\perp \Big( \frac1{2\pi} \Big(1- \sum_{i,j} \varphi_{i,j}^\eps(x) \Big) \int_{\Omega^\eps} \ln|x-y| f(y)\, d y\Big)
\]
which is divergence free, vanishing close to the boundary. Hence it verifies \eqref{eq.veps}, but in this case, we could only have
\[
 \| K_{\R^2}[f]-\tilde v^\eps[f] \|_{L^2(\Omega^\eps)} \leq C \|f\|_{L^1\cap L^\infty} \Bigl( \| \nabla \varphi^\eps \|_{L^2} + \| \varphi^\eps \|_{L^2} \Bigl) \sqrt{\sum_{i,j} 1},
\]
which tends to zero only for large distances (see the previous remark).

The use of the Biot Savart law in the exterior of one inclusion allowed us to get the smallness of the stream function $\psi^\eps$ on the support of $\nabla \varphi_{i,j}^\eps$. Nevertheless, there is another classical way to construct a compatible correction, as it used in \cite{Allaire90a,Allaire90b,Tartar80} and more recently in \cite{LaMa}. The general idea is to utilize a direct cut-off function on the velocity and a
correction to restore the divergence-free condition, based on Bogovski{\u\i} operator. Namely, let us consider the case $d_{\varepsilon}\geq C_{0}\varepsilon$ for some $C_{0}>0$, and we use $\varphi^\varepsilon$ independent of $d_{\varepsilon}$ (see Case 1, just above). Then we define:
\[
\hat v^\varepsilon[f](x):= \Big(1- \sum_{i,j} \varphi_{i,j}^\eps(x) \Big) K_{\R^2}[f] + h^\varepsilon
\]
where $h^\varepsilon$ is a solution of
\[
\div h^\varepsilon = K_{\R^2}[f] \cdot \sum_{i,j}\nabla \varphi_{i,j}^\eps(x) . 
\]
By a dilatation argument and standard estimates for such a problem, we can prove the existence of $h^\varepsilon$ which satisfies the same estimate than $ \varphi_{i,j}^\eps$, namely Lemma 2.1 in \cite{LaMa} gives
\[
\| h^\varepsilon \|_{L^2(\Omega^\eps)} \leq C \|f\|_{L^1\cap L^\infty} \varepsilon \| \nabla \varphi^\eps \|_{L^2} \sqrt{\sum_{i,j} 1}.
\]
Therefore, without using Biot Savart law in the exterior of one inclusion, we have in a much simpler way
\[
 \| K_{\R^2}[f]-\hat v^\eps[f] \|_{L^2(\Omega^\eps)} \leq C \|f\|_{L^1\cap L^\infty} \Bigl( \eps\| \nabla \varphi^\eps \|_{L^2} + \| \varphi^\eps \|_{L^2} \Bigl) \sqrt{\sum_{i,j} 1},
\]
which gives the result when $d_{\varepsilon}\geq C_{0}\varepsilon$. In particular, this correction is sufficient for the case where the inclusions cover the square. This argument cannot be used when the distance is smaller than the size of the inclusions. For distances $d^\eps \ll 
\eps$, we would need to understand the behavior of the estimates of $h^\varepsilon$ on domains of the form $ (-1,1)^2\setminus \rho \K$ 
as $\rho\to 1^-$, which is a complicated question. This explains why \cite{Allaire90b,LaMa} need to assume that size is still smaller than the inter-hole distance (see \cite[Equation (3.4.2)]{Allaire90b} and \cite[Remark 2.1]{LaMa}).

To conclude, the use of the Biot Savart law in the exterior of one inclusion is necessary for distances smaller than the size of the inclusions, as in the case of inclusions distributed on the segment. As the estimates of $w^{k,\varepsilon}_{i,j}$ are the same for both cases, we have chosen here to keep the correction $v^\varepsilon$, even in the case where inclusions cover the square.
\end{remark}

\section{Impermeability for small distance}\label{sec:impermeability}

The goal of this section is to prove Proposition~\ref{prop:key2}.
We begin with the easier situation which is the case of inclusions distributed on the segment \eqref{Omega1eps}.

\subsection{Inclusions distributed on the segment} 
In this subsection we consider $\Omega^\eps_{1}$ in the case $d_{\eps} \leq \eps$ and a bounded sequence in $L^2_{\loc}(\R^2)$, of divergence free vector fields $v^\eps$ which are tangent to the boundary $\partial\Omega^\eps_{1}$. Let us fix $p\in C^\infty_{c} (\overline{\Omega_{1}})$.

For any $s\in (0,\rho_{0}]$, we connect all the inclusions
by taking
\[
\mathcal{C}^\varepsilon(s):=\Big(\bigcup_{i=1}^{N_{\eps}} \K_{i,1}^\varepsilon\Big) \bigcup \Big([x_{1,1},x_{N_{\eps},1}]\times [-\tfrac\eps2 s,\tfrac\eps2 s]) \Big),
\]
where $x_{i,j}$ is the horizontal coordinate of $z_{i,j}$. 

Due to the assumption (H2) on $\K$, we have
\begin{itemize}
\item $\mathcal{C}^\varepsilon(s)$ is a simply connected compact subset of $\R^2$;
\item $\mathcal{C}^\varepsilon(s)\setminus \Big(\bigcup_{i=1}^{N_{\eps}} \K_{i,1}^\varepsilon\Big)=\Omega^\varepsilon_{1} \cap \mathcal{C}^\varepsilon(s)$ has $N_{\eps}-1$ connected components whose Lebesgue measure can be 
estimated as follows:
\[
\mathcal{A}^\varepsilon(s) := {\rm meas} \Big|\Omega^\varepsilon \cap \mathcal{C}^\varepsilon(s)\Big| \leq (N_{\eps}-1) \Big( \Big(\frac{\varepsilon}2\Big)^2 4 \int_{0}^{s} \rho_{1} r^{1+\gamma}\, dr + \eps s d_{\eps}\Big) 
\leq C ( \varepsilon s^{\gamma +2}+ d_{\eps} s),
\]
with $C$ independent of $\varepsilon$ (where we have used that $N_{\eps}\leq 2/\eps$).
\end{itemize}
See the left hand side picture of Figure~\ref{fig.area} (page \pageref{fig.area}) to understand the space between the inclusions. 

Unfortunately, we cannot state that the point $(1,0)$ belongs to $\mathcal{C}^\varepsilon(s)$. For this reason, we introduce a cutoff function
\begin{equation}\label{chi}
\chi \in C^\infty (\R^+),\quad \chi (r)=0 \text{ if }r\geq 2 \text{ and } \chi (r)=1 \text{ if }r\leq 1
\end{equation}
and for any $\eta >0$ we decompose $p$ as
$$
p(x) = \Big(1-\chi\Big(\frac{|x-(1,0)|}\eta\Big)\Big)p(x) +\chi\Big(\frac{|x-(1,0)|}\eta\Big)p(x)=:p_{1,\eta}(x)+p_{2,\eta}(x).
$$

Before starting the actual proof, let us note that the main idea of the following analysis is to observe that 
 $$
 \int_{\Omega^\eps_{1}\cap \Omega_{1}} v^\eps\cdot \nabla p_{1,\eta}\, dx
\approx \int_{\R^2 \setminus \mathcal{C}^\varepsilon(s) } v^\eps\cdot \nabla p_{1,\eta}\, dx
= \int_{ \partial \mathcal{C}^\varepsilon(s)\setminus \partial \Omega^\varepsilon_{1} } p_{1,\eta} v^\eps \cdot n\, dx
\quad \forall s.
$$
Next we will integrate this equality on $s\in [0,s_{\varepsilon}]$ for some $s_{\varepsilon}$ chosen later, and we estimate it by $\|p_{1,\eta}\|_{L^\infty}\| v^\eps \|_{L^2} (\mathcal{A}^\varepsilon(s_{\varepsilon}))^{1/2}$. We will conclude by finding the best $s_{\varepsilon}$ (depending on $\gamma$) such that $ \int_{\Omega^\eps_{1}\cap \Omega_{1}} v^\eps\cdot \nabla p_{1,\eta}\, dx$ tends to zero as $\varepsilon\to 0$.

More precisely, let us fix $\eta>0$ and prove the convergence with $p_{1,\eta}$. Thanks to (H2), we state that for $\varepsilon$ small enough (namely, $\varepsilon<\eta$) we have 
\begin{itemize}
\item $([0,1-\eta]\times\{ 0\}) \subset \mathcal{C}^\varepsilon(s)$ for any $s$.
\end{itemize}
As $p_{1,\eta}\equiv 0$ in $B((1,0),\eta)$, we deduce that $p_{1,\eta}$ is smooth in $\R^2 \setminus \mathcal{C}^\varepsilon(s)$ and may have a jump across $([0,1-\eta]\times\{ 0\})$. This remark is crucial for the following computation (and this explains the introduction of the cutoff). Using that $v^\eps$ is tangent to $\partial \Omega^\varepsilon$ and that $p_{1,\eta}$ is smooth on $\R^2 \setminus \mathcal{C}^\varepsilon(s)$, we compute for any $s\in (0,s_{\eps})$, with $s_{\eps}$ to be determined later such that $s_{\eps}\leq \rho_{0}$:
\begin{eqnarray*}
 \int_{\Omega^\eps_{1}\cap \Omega_{1}} v^\eps\cdot \nabla p_{1,\eta}\, dx
&=& \int_{\R^2 \setminus \mathcal{C}^\varepsilon(s) } v^\eps\cdot \nabla p_{1,\eta}\, dx + \int_{\Omega^\varepsilon_{1} \cap \Omega_{1} \cap \mathcal{C}^\varepsilon(s)} v^\eps\cdot \nabla p_{1,\eta}\, dx\\
&=& \int_{ \partial \mathcal{C}^\varepsilon(s)\setminus \partial \Omega^\varepsilon_{1} } p_{1,\eta} v^\eps \cdot n\, dx + \int_{\Omega^\varepsilon_{1} \cap \Omega_{1} \cap \mathcal{C}^\varepsilon(s)} v^\eps\cdot \nabla p_{1,\eta}\, dx 
\end{eqnarray*}
where $n$ is equal to $\pm e_{2}$. Moreover, $ \partial \mathcal{C}^\varepsilon(s)\setminus \partial \Omega^\varepsilon_{1}$ has $2(N_{\eps}-1)$ connected components, which are horizontal segments linking $\K^\varepsilon_{i,1}$ and $\K^\varepsilon_{i+1,1}$, and included in the lines $([0,1]\times \{\frac\varepsilon2 s\})$ and $([0,1]\times \{-\frac\varepsilon2 s\})$. Now, we integrate the above equality for $s\in (0,s_{\eps})$:
\begin{equation}\label{main compute}
s_{\eps} \int_{\Omega_{1}^\eps\cap \Omega_{1}} v^\eps\cdot \nabla p_{1,\eta}\, dx
= \int_{0}^{s_{\eps}} \int_{ \partial \mathcal{C}^\varepsilon(s)\setminus \partial \Omega_{1}^\varepsilon } p_{1,\eta} v^\eps \cdot n\, dxds + \int_{0}^{s_{\eps}} \int_{\Omega_{1}^\varepsilon \cap \Omega_{1} \cap \mathcal{C}^\varepsilon(s)} v^\eps\cdot \nabla p_{1,\eta}\, dxds.
\end{equation}
Changing variable $s'=\frac{\varepsilon s}2$, the first right hand side term can be estimated as follows:
\begin{eqnarray*}
s_{\eps}^{-1} \Big| \int_{0}^{s_{\eps}} \int_{ \partial \mathcal{C}^\varepsilon(s)\setminus \partial \Omega_{1}^\varepsilon } p_{1,\eta} v^\eps \cdot n\, dxds \Big|&=&2(\varepsilon s_{\eps})^{-1}
\Big| \int_{0}^{\varepsilon s_{\eps}/2}\int_{ \partial \mathcal{C}^\varepsilon(2s'/\varepsilon) \setminus \partial \Omega_{1}^\varepsilon } p_{1,\eta} v^\eps \cdot n\, dxds' \Big|\\
 &\leq & 2 (\varepsilon s_{\eps})^{-1} \int_{\Omega_{1}^\varepsilon \cap \mathcal{C}^\varepsilon(s_{\eps}) } |p_{1,\eta}| |v^\eps| \, dx \\
 &\leq &2 (\varepsilon s_{\eps})^{-1} \| p_{1,\eta} \|_{L^\infty(\Omega_{1})} \| v^\eps \|_{L^2(\Omega_{1}^\eps \cap \ \supp p)} \Big(\mathcal{A}^\varepsilon(s_{\eps})\Big)^{1/2}.
\end{eqnarray*}
Now we choose $s_{\eps}$ such that $\frac{\mathcal{A}^\varepsilon(s_{\eps})}{(\eps s_{\eps})^2}\leq C( \frac{s_{\eps}^\gamma}\eps + \frac{d_{\eps}}{\eps^{2}s_{\eps}})$ is minimal, i.e. $s_{\eps}= (\frac{d_{\eps}}{\gamma \eps})^{1/(\gamma+1)}$ for $\gamma<\infty$ (which is smaller than $\rho_{0}$ for $\eps$ small enough) and $s_{\eps} =\rho_{0}$ for $\gamma=\infty$. 

In both case, we have
\[
s_{\eps}^{-1} \Big| \int_{0}^{s_{\eps}} \int_{ \partial \mathcal{C}^\varepsilon(s)\setminus \partial \Omega_{1}^\varepsilon } p_{1,\eta} v^\eps \cdot n\, dxds \Big| \leq C \| p_{1,\eta} \|_{L^\infty(\Omega_{1})} \| v^\eps \|_{L^2(\Omega_{1}^\eps \cap \ \supp p)} \frac1{\sqrt{\eps}}\Big(\frac{d_{\eps}}{\eps}\Big)^{\frac{\gamma}{2(\gamma+1)}}.
\]

For the second right hand side term of \eqref{main compute}, we have
\begin{eqnarray*}
s_{\eps}^{-1}\Big| \int_{0}^{s_{\eps}} \int_{\Omega_{1}^\varepsilon \cap \Omega_{1}\cap \mathcal{C}^\varepsilon(s)} v^\eps\cdot \nabla p_{1,\eta}\, dxds \Big|
&\leq& \int_{\Omega_{1}^\varepsilon\cap \Omega_{1} \cap \mathcal{C}^\varepsilon(s_{\eps})} | v^\eps| |\nabla p_{1,\eta} | \, dx\\
&\leq& \| \nabla p_{1,\eta} \|_{L^\infty(\Omega_{1})} \| v^\eps \|_{L^2(\Omega_{1}^\eps \cap \ \supp p)} (\varepsilon s_{\eps} )(\varepsilon s_{\eps} )^{-1}\Big(\mathcal{A}^\varepsilon(s_{\eps})\Big)^{1/2}\\
&\leq& C\| p \|_{W^{1,\infty}(\Omega_{1})}\| v^\eps \|_{L^2(\Omega_{1}^\eps \cap \ \supp p)} d_{\eps}^{\frac1{\gamma+1}} \eps^{\frac\gamma{\gamma+1}} \frac1{\sqrt{\eps}}\Big(\frac{d_{\eps}}{\eps}\Big)^{\frac{\gamma}{2(\gamma+1)}} .
\end{eqnarray*}

Bringing together these two estimates with \eqref{main compute}, we have
\begin{eqnarray*}
\Big| \int_{\Omega_{1}^\eps \cap \Omega_{1}} v^\eps\cdot \nabla p_{1,\eta}\, dx \Big|
&\leq& C\| p \|_{W^{1,\infty}(\Omega_{1})}\| v^\eps \|_{L^2(\Omega_{1}^\eps \cap \ \supp p)} \frac1{\sqrt{\eps}}\Big(\frac{d_{\eps}}{\eps}\Big)^{\frac{\gamma}{2(\gamma+1)}} \\
&\leq &C\| p \|_{W^{1,\infty}(\Omega_{1})}\| v^\eps \|_{L^2(\Omega_{1}^\eps \cap \ \supp p)} \Bigg(\frac{d_{\eps}}{\eps^{2+\frac1\gamma}}\Bigg)^{\frac{\gamma}{2(\gamma+1)}},
\end{eqnarray*}
which tends to zero (under the assumptions of Proposition~\ref{prop:key2}).

We conclude now with $p$. As $v^\eps$ is bounded in $L^2(\Omega_{1}\cap \supp p)$, we can extract a subsequence such that its converges weakly to $v$. So, for every $\eta>0$, 
\[
\int_{\Omega_{1}} v^\eps \cdot \nabla p_{2,\eta} \, dx \to \int_{\Omega_{1}} v \cdot \nabla p_{2,\eta} \, dx \quad \text{as } \eps\to 0.
\]
Moreover, as $\nabla p_{2,\eta}$ converges weakly to zero in $L^2$ when $\eta\to 0$, we conclude that for any $\delta >0$, there exists $\eta$ such that
\[
\Big|\int_{\Omega_{1}} v \cdot \nabla p_{2,\eta} \, dx\Big | \leq \delta/3
\]
and for this $\eta$, there exists $\eps_{\eta}$ such that for any $\eps \in (0,\eps_{\eta}]$:
\begin{align*}
 \Big| \int_{\Omega_{1}} v^\eps \cdot \nabla p \, dx \Big| &\leq \Big| \int_{\Omega_{1}} v^\eps \cdot \nabla p_{1,\eta} \, dx \Big| + \Big| \int_{\Omega_{1}} v^\eps \cdot \nabla p_{2,\eta} \, dx \Big| \\
 & \leq \frac{\delta}3 + \Big| \int_{\Omega_{1}} v \cdot \nabla p_{2,\eta} \, dx \Big| + \frac{\delta}3 \leq \delta.
\end{align*}
Therefore, we have proved that for any sequence $\eps_{n}\to 0$, we can extract a subsequence such that $ \int_{\Omega_{1}} v^{\eps_{\psi(n)}} \cdot \nabla p \, dx$ tends to zero. By uniqueness of this limit, we deduce that the above convergence holds for the full sequence, without extraction. 

This ends the proof of the first point of Proposition~\ref{prop:key2}.

\subsection{Inclusions distributed on the square}\label{imperm:square}

We consider $\Omega_{2}^\eps$ \eqref{Omega2eps} in the case $d_{\eps} \leq \eps$ and a bounded sequence in $L^2_{\loc}(\R^2)$, of divergence free vector fields $v^\eps$ which are tangent to the boundary $\partial\Omega^\eps_{2}$. Let us fix $p\in C^\infty_{c} (\overline{\Omega_{2}})$ which is extended by zero in $(0,1)^2$.

 We use the cutoff \eqref{chi} in the neighborhood of the four corners: for any $\eta >0$ fixed, we decompose $p$ as
\begin{equation*}
\begin{split}
p(x) = &\Big(1-\chi\Big(\frac{|x-(1,0)|}\eta\Big)-\chi\Big(\frac{|x|}\eta\Big)-\chi\Big(\frac{|x-(0,1)|}\eta\Big)-\chi\Big(\frac{|x-(1,1)|}\eta\Big)\Big)p(x) \\
&+\Big(\chi\Big(\frac{|x-(1,0)|}\eta\Big)+\chi\Big(\frac{|x|}\eta\Big)+\chi\Big(\frac{|x-(0,1)|}\eta\Big)+\chi\Big(\frac{|x-(1,1)|}\eta\Big)\Big)p(x)\\ 
=:&p_{1,\eta}(x)+p_{2,\eta}(x),
\end{split}
\end{equation*}
and the core of the analysis is to prove that
\begin{equation*}
\int_{\Omega_{2}} v^\eps \cdot \nabla p_{1,\eta} \, dx \to 0\quad \text{as } \eps\to 0.
\end{equation*}
Up to using other cutoffs, we can assume that $p_{1,\eta}$ is equal to zero in the neighborhood of three sides of the square. For instance, let us assume that $p_{1,\eta}$ vanishes except around the bottom side $(0,1) \times \{ 0 \}$. 

As $\partial \K$ is $C^1$, there exist $\gamma>0$, $\rho_{0},\rho_{1}>0$ such that
\[
\Big([-(1-\rho_{1} |s|^{1+\gamma}), -(1-\rho_{1} \rho_{0}^{1+\gamma})]\cup [(1-\rho_{1} \rho_{0}^{1+\gamma}), (1-\rho_{1} |s|^{1+\gamma})]\Big) \times\{s\}\subset \K
\]
for all $s\in [-\rho_{0},\rho_{0}]$. For simplicity, let assume that the entire segment $[-(1-\rho_{1} |s|^{1+\gamma}),(1-\rho_{1} |s|^{1+\gamma})]\times\{s\}\subset \K$ for all $s\in [-\rho_{0},\rho_{0}]$ (see the discussion in the final section).

The main idea is to see the side $[0,1]\times \{0\}$ as the limit of $2 N'+1$ lines of inclusions where we set
\begin{equation}\label{Nprime}
N'=\Big[ \Big(\frac{d_{\eps}}{\eps}\Big)^{\frac\gamma{2(\gamma+1)}} \frac1\eps \Big], 
\end{equation}
and to use the area between $\K_{i,j}^\varepsilon$ and $\K_{i+1,j}^\varepsilon$ in each lines $j=1,\dots,2N'+1$.

Let us first note that the case $N'=0$, i.e. $\frac{d_{\eps}}{\eps^{2+\frac1\gamma}}<\eps$ is already covered by the previous analysis: the fluid cannot penetrated though one line of obstacles. So, for the rest of this section, we consider the case $N'\geq 1$. We also remark that $2N'+1 \ll N_{\eps} = \mathcal{O}(1/\eps)$ because we assume in this subsection that $\frac{d_{\eps}}\eps\to 0$. Hence we have for $\varepsilon$ small enough that
\begin{itemize}
\item $\K_{1,j}^\varepsilon \subset B(0,\eta)$ and $\K_{N_{\eps},j}^\varepsilon \subset B((1,0),\eta)$ for all $j=1,\dots, 2N'+1$.
\end{itemize}

For any $s\in [-\rho_{0},\rho_{0}]$ and $k=1,\dots, N'$, we connect the inclusions between the lines $N'+1-k$ and $N'+1+k$:
\[
\mathcal{C}^\varepsilon_{k}(s):=\Big(\bigcup_{i=1}^{N_{\eps}}\bigcup_{j=N'+1-k}^{N'+1+k} \K_{i,j}^\varepsilon\Big) \bigcup \Big([x_{1,1},x_{N_{\eps},1}]\times [y_{1,N'+1-k}-\tfrac\eps2 s,y_{1,N'+1+k}+\tfrac\eps2 s]) \Big).
\]
where $y_{i,j}$ is the vertical coordinate of $z_{i,j}$. Hence, we have for any $k=1,\dots, N'$:
\begin{itemize}
\item for all $s\in [-\rho_{0},\rho_{0}]$, $\mathcal{C}_{k}^\varepsilon(s)$ is a simply connected compact subset of $\R^2$;
\item for all $s\in (0,\rho_{0}]$, $\Bigl(\mathcal{C}_{k}^\varepsilon(s)\setminus \mathcal{C}_{k}^\varepsilon(-s)\Big) \setminus \Big(B(0,\eta)\cup B((1,0),\eta)\Big)$ has less than $2(N_{\eps}-1)$ connected components whose Lebesgue measure can be estimated as follows:
\begin{equation*}
\begin{split}
\mathcal{A}_{k}^\varepsilon(s) &:= {\rm meas} \Big|\Bigl(\mathcal{C}_{k}^\varepsilon(s)\setminus \mathcal{C}_{k}^\varepsilon(-s)\Big) \setminus \Big(B(0,\eta)\cup B((1,0),\eta)\Big)\Big|\\
& \leq 2 (N_{\eps}-1) \Big( \Big(\frac\varepsilon4\Big)^2 4 \int_{0}^{s} \rho_{1} r^{1+\gamma}\, dr + \eps sd_{\eps}\Big) 
\leq C ( \varepsilon s^{\gamma +2}+ d_{\eps} s),
\end{split}
\end{equation*}
with $C$ independent of $\varepsilon$.
\end{itemize}
As in the case of one line of obstacles, we can show that the best choice is $s_{\eps} =(\frac{d_{\eps}}{\gamma\eps})^{1/(\gamma+1)}$ for $\gamma<\infty$ (which belongs to $(0,\rho_{0})$ for $\eps$ small enough) and $s_{\eps}=\rho_{0}$ for $\gamma=\infty$. Hence we have
\begin{equation}\label{aire2}
\mathcal{A}^\varepsilon_{k}(s_{\eps}) \leq 2C d_{\eps}s_\eps.
\end{equation}
Here too, it is important for the following argument that the jump of the test function can occur only in $\mathcal{C}_{k}^\varepsilon(s)$, so let us note that
\begin{itemize}
 \item $\tilde p_{1,\eta}:=p_{1,\eta}(\cdot - (0, y_{1,N'+1}))$ is smooth on $\R^2\setminus \mathcal{C}_{k}^\varepsilon(s_{\eps})$ for all $s\in [-\rho_{0},\rho_{0}]$ and $k=1,\dots, N'$.
\end{itemize}

Now, we use that $v^\eps$ is tangent to $\partial \Omega^\varepsilon_{2}$ to compute for any $s\in (-s_{\eps},s_{\eps})$ and $k=1,\dots,N'$:
\begin{eqnarray*}
 \int_{\Omega_2^\eps} v^\eps\cdot \widehat{\nabla \tilde p_{1,\eta}}\, dx
&=& \int_{\R^2 \setminus \mathcal{C}_{k}^\varepsilon(s) } v^\eps\cdot \nabla \tilde p_{1,\eta}\, dx + \int_{\Omega_2^\varepsilon \cap \mathcal{C}_{k}^\varepsilon(s)} v^\eps\cdot \widehat{\nabla \tilde p_{1,\eta}}\, dx\\
&=& \int_{ \partial \mathcal{C}_{k}^\varepsilon(s)\setminus \partial \Omega_2^\varepsilon } \tilde p_{1,\eta} v^\eps \cdot n\, dx + \int_{\Omega_2^\varepsilon \cap \mathcal{C}_{k}^\varepsilon(s)} v^\eps\cdot \widehat{\nabla \tilde p_{1,\eta}}\, dx 
\end{eqnarray*}
where $n$ is equal to $\pm e_{2}$ and where $\widehat{\nabla \tilde p_{1,\eta}}$ is the extension of $\nabla \tilde p_{1,\eta}$ by zero for $y\geq y_{1,N'+1}$. Moreover, $ (\partial \mathcal{C}_{k}^\varepsilon(s)\setminus \partial \Omega_2^\varepsilon)\cap {\rm supp\ }\tilde p_{1,\eta}$ has less than $2(N_{\eps}-1)$ connected components, which are horizontal segments linking $\K^\varepsilon_{i,j}$ and $\K^\varepsilon_{i+1,j}$ for $j=N'+1\pm k$, and included in the lines $([0,1]\times \{y_{1,N'+1+k}+\varepsilon s\})$ and $([0,1]\times \{y_{1,N'+1-k}-\varepsilon s\})$. Now, we integrate the above equality for $s\in (-s_{\eps},s_{\eps})$ and sum on $k=1,\dots, N'$
\begin{equation}\label{main compute2}
\begin{split}
N'(2s_{\eps}) \int_{\Omega_2^\eps} v^\eps\cdot \widehat{\nabla \tilde p_{1,\eta}}\, dx
=& \sum_{k=1}^{N'} \int_{-s_{\eps}}^{s_{\eps}} \int_{\partial \mathcal{C}_{k}^\varepsilon(s)\setminus \partial \Omega_2^\varepsilon } \tilde p_{1,\eta} v^\eps \cdot n\, dxds \\
&+ \sum_{k=1}^{N'} \int_{-s_{\eps}}^{s_{\eps}} \int_{\Omega_2^\varepsilon \cap \mathcal{C}_{k}^\varepsilon(s)} v^\eps\cdot \widehat{\nabla \tilde p_{1,\eta}}\, dxds.
\end{split}
\end{equation}
Changing variable $s'=\frac{\varepsilon s}2$, the first right hand side term can be estimated as follows:
\begin{eqnarray*}
 \Big| \int_{-s_{\eps}}^{s_{\eps}} \int_{ \partial \mathcal{C}_{k}^\varepsilon(s)\setminus \partial \Omega_2^\varepsilon } \tilde p_{1,\eta} v^\eps \cdot n\, dxds \Big|&=
& 2\eps ^{-1} \Big| \int_{-\eps s_{\eps}/2}^{\eps s_{\eps}/2}\int_{ \partial \mathcal{C}_{k}^\varepsilon(2s'/\varepsilon) \setminus \partial \Omega_2^\varepsilon } \tilde p_{1,\eta} v^\eps \cdot n\, dxds' \Big|\\
 &\leq & 2\eps^{-1} \int_{ \mathcal{C}_{k}^\varepsilon(s_{\eps})\setminus \mathcal{C}_{k}^\varepsilon(-s_{\eps}) } |\tilde p_{1,\eta}| |v^\eps| \, dx .
 \end{eqnarray*}
As $(\mathcal{C}_{k}^\varepsilon(s_{\eps})\setminus \mathcal{C}_{k}^\varepsilon(-s_{\eps}))\cap (\mathcal{C}_{k'}^\varepsilon(s_{\eps})\setminus \mathcal{C}_{k'}^\varepsilon(-s_{\eps}))\cap ({\rm supp\ }\tilde p_{1,\eta})=\emptyset$ if $k'\neq k$, the sum becomes:
\begin{align*}
(N's_{\eps})^{-1} \Big| \sum_{k=1}^{N'} \int_{-s_{\eps}}^{s_{\eps}} \int_{ \partial \mathcal{C}_{k}^\varepsilon(s)\setminus \partial \Omega_2^\varepsilon } \tilde p_{1,\eta} v^\eps \cdot n\, dxds \Big| 
 &\leq (N's_{\eps} \eps)^{-1} \| p_{1,\eta} \|_{L^\infty(\Omega_2)} \| v^\eps \|_{L^2(\Omega_2^\varepsilon\cap \supp \tilde p_{1,\eta})} \Big(\sum_{k=1}^{N'} \mathcal{A}_{k}^\varepsilon(s_{\eps})\Big)^{1/2}\\
 &\leq C\Big(\frac{d_{\eps}}{N's_{\eps} \eps^2} \Big)^{1/2} \| p_{1,\eta} \|_{L^\infty(\Omega_2)} \| v^\eps \|_{L^2(\Omega_2^\varepsilon\cap \supp \tilde p_{1,\eta})}\\
 &\leq C\sqrt{2}\Big(\frac{d_{\eps}}{(N'+1)s_{\eps} \eps^2} \Big)^{1/2} \| p_{1,\eta} \|_{L^\infty(\Omega_2)} \| v^\eps \|_{L^2(\Omega_2^\varepsilon\cap \supp \tilde p_{1,\eta})}\\
 &\leq C\Big(\frac{d_\eps}{\eps}\Big)^{\frac\gamma{4(\gamma+1)}} \| p_{1,\eta} \|_{L^\infty(\Omega_2)} \| v^\eps \|_{L^2(\Omega_2^\varepsilon\cap \supp \tilde p_{1,\eta})}
\end{align*}
where we have used \eqref{aire2} and \eqref{Nprime}. We can remark in the previous computation that we introduce cutoff functions close to the corners of the square to avoid integrals between $\K_{1,j}^\eps$ and $\K_{1,j+1}^\eps$.

For the second right hand side term of \eqref{main compute2}, we state
\begin{equation*}\begin{split}
\Big| \sum_{k=1}^{N'} \int_{-s_{\eps}}^{s_{\eps}} \int_{\Omega_2^\varepsilon \cap \mathcal{C}_{k}^\varepsilon(s)} v^\eps\cdot \widehat{\nabla \tilde p_{1,\eta}}\, dxds \Big|
\leq& N' (2s_{\eps}) \int_{\Omega_2^\varepsilon \cap \mathcal{C}_{N'}^\varepsilon(s_{\eps})} | v^\eps| |\widehat{\nabla \tilde p_{1,\eta}} | \, dx\\
\leq N' (2s_{\eps})& \| \nabla p_{1,\eta} \|_{L^\infty(\Omega_2)}\| v^\eps \|_{L^2(\Omega_2^\varepsilon\cap \supp \tilde p_{1,\eta})} \Big( C(\varepsilon+d_{\eps})N' \Big)^{1/2}\\
\leq N' (2s_{\eps})& C \| p_{1,\eta} \|_{W^{1,\infty}(\Omega_2)}\| v^\eps \|_{L^2(\Omega_2^\varepsilon\cap \supp \tilde p_{1,\eta})} \Big(\frac{d_\eps}{\eps}\Big)^{\frac\gamma{4(\gamma+1)}}.
\end{split}\end{equation*}

Bringing together these two estimates with \eqref{main compute2}, we have
\[
\Big| \int_{\Omega_2^\eps} v^\eps\cdot \widehat{\nabla \tilde p_{1,\eta}}\, dx \Big|
\leq C \| p_{1,\eta} \|_{W^{1,\infty}(\Omega_2)}\| v^\eps \|_{L^2(\Omega_2^\varepsilon\cap \supp \tilde p_{1,\eta})} \Big(\frac{d_\eps}{\eps}\Big)^{\frac\gamma{4(\gamma+1)}},
\]
which goes to zero as $\eps \to 0$ (under the assumption of Proposition~\ref{prop:key2}).

Hence, we deduce that
\[
\Big|\int_{\Omega_{2}} v^\eps \cdot \nabla p_{1,\eta} \, dx \Big|
\leq \Big| \int_{\Omega_{2}^\eps} v^\eps \cdot \widehat{\nabla \tilde p_{1,\eta}} \, dx \Big|
+
\| v^\eps \|_{L^2( \supp p_{1,\eta} \cup \supp \tilde p_{1,\eta})} \| \widehat{\nabla p_{1,\eta}}(\cdot)-\widehat{\nabla p_{1,\eta}}( \cdot - (0, y_{1,N'+1})) \|_{L^2}
\]
tends to zero when $\eps \to 0$, because $y_{1,N'+1}\leq 2\eps (N'+1)\to 0$. In the previous estimate, $\widehat{\nabla p_{1,\eta}}$ is equal to $\nabla p_{1,\eta}$ for $y<0$ and to zero for $y\geq 0$, which belongs to $L^2$.

Finally, we conclude with $p=p_{1,\eta}(x)+p_{2,\eta}$ exactly in the way than the case where the obstacles are distributed on the segment: for any $\delta >0$, there exists $\eta$ and $\eps_{\eta}$ such that for any $\eps \in (0,\eps_{\eta}]$:
\begin{align*}
 \Big| \int_{\Omega_{2}} v^\eps \cdot \nabla p \, dx \Big| &\leq \Big| \int_{\Omega_{2}} v^\eps \cdot \nabla p_{1,\eta} \, dx \Big| + \Big| \int_{\Omega_{2}} v^\eps \cdot \nabla p_{2,\eta} \, dx \Big| \\
 & \leq \frac{\delta}3 + \Big| \int_{\Omega_{2}} v \cdot \nabla p_{2,\eta} \, dx \Big| + \frac{\delta}3 \leq \delta
\end{align*}
This ends the proof of Proposition~\ref{prop:key2}.

\section{Final remarks and comments}\label{sect final}

\subsection{Weaker assumptions on the shape}\label{sect shape}

In this subsection, we discuss how to decrease the assumption (H2).

The goal of (H2) is to construct a good cutoff function supported between the inclusions (see Subsection~\ref{subsec:cutoff}) for the permeability result, and to compute the area between two inclusions (see Section~\ref{sec:impermeability}) for the impermeability result.

\begin{figure}[h!t]
\begin{center}
\includegraphics[height=5cm]{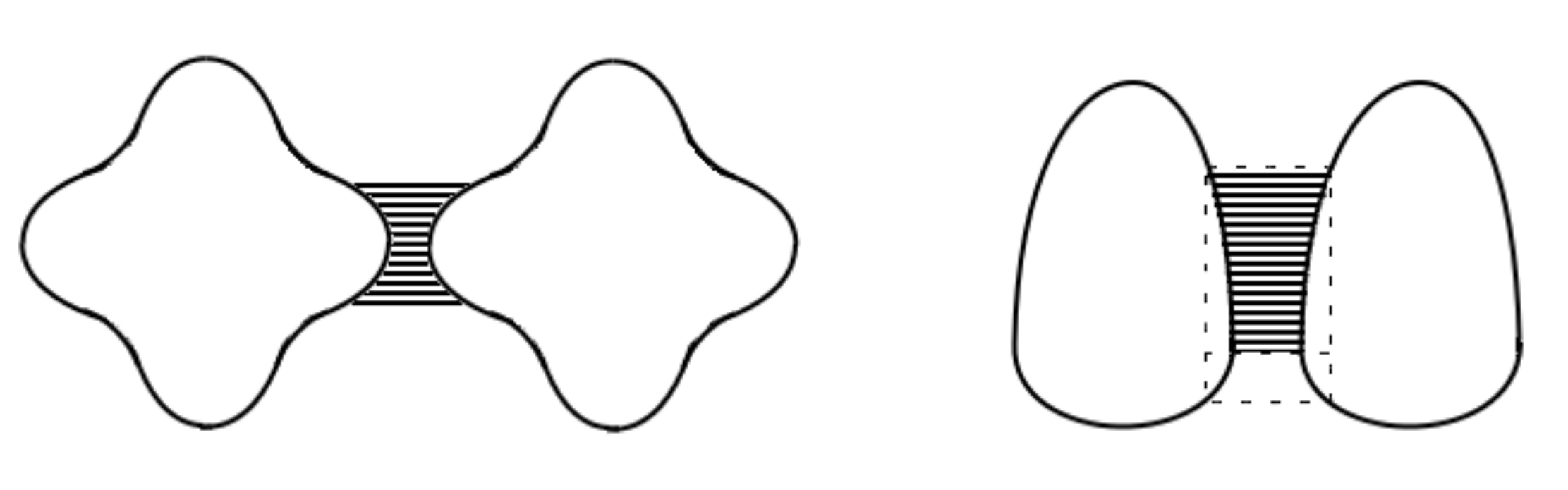}
\caption{Area between two inclusions.}\label{fig.area}
\end{center}
\end{figure}

For both arguments, it is not necessary to assume:
\begin{itemize}
\item the symmetry assumption on $\partial\K$ around $(1,0)$ (it was written only for clarity). For example, let us consider the case where $[-(1-\rho_{1} |s|^{1+\gamma_{1}}),(1-\rho_{1} |s|^{1+\gamma_{1}})]\times\{s\}\subset \K$ for all $s\in [0,\rho_{0}]$ and $[-(1-\rho_{2} |s|^{1+\gamma_{2}}),(1-\rho_{2} |s|^{1+\gamma_{2}})]\times\{s\}\subset \K$ for all $s\in [-\rho_{0},0]$ with $\gamma_{1}>\gamma_{2}$ (see the right hand side picture of Figure~\ref{fig.area}). If we perform the proof of Proposition~\ref{prop:key2} with the area between the lines $y=\frac\eps2 s_{\eps}$ and $y= -\frac\eps2 s_{\eps}$, we would obtain the impermeability as if $\gamma=\gamma_{2}=\min (\gamma_{1},\gamma_{2})$ which is not optimal. At the opposite, if we perform the proof with the area between the lines $y=\frac\eps2 s_{\eps}$ and $y= 0$, then we get that the asymptotic behavior depends on the limit of $d_{\eps}/\eps^{2+\frac{1}{\gamma_{1}}}$.

\item that the segment $[(-1,0),(1,0)]\subset \K$. This assumption simplified the construction of $\mathcal{C}^\eps(s)$ with is the connexion of $\K_{i,1}$ (see Section~\ref{sec:impermeability}). Without this assumption, we should define 
$\mathcal{C}^\eps(s):=\cup\K_{i,1} \cup A_{i}$, 
where $A_{i}$ is the area between $\K_{i,1}$ and $\K_{i+1,1}$.

\item that $(\pm 1,0),(0,\pm 1)\in \partial\K$. Without this assumption, we have to change sightly \eqref{domain2} such that the distance between $\K_{i,j}$ and $\K_{i+1,j}$ (and between $\K_{i,j}$ and $\K_{i,j+1}$) is $d_{\eps}$. Next, we should consider the coefficient $\gamma$ around the closest points, with respect to an axis which can be not horizontal. In the same spirit, we can also distribute the inclusions on any Jordan arc or on any compact set with no empty interior.

\item that the inclusions have the same shape (see \eqref{domain1}). For instance, we can consider a finite number of possible shapes $(\K^1,\dots, \K^M)$, and we distribute the inclusions on the square:
\[\K_{i,j}^{\eps}:= z_{i,j}^{\varepsilon} + \tfrac\eps 2 \K^{p_{i,j}} \text{ with }p_{i,j}\in \{1,\dots, M\}.\]
In this case, considering uniform constants to estimate every conformal mappings $\Tc^1,\dots,\Tc^M$, we can adapt the permeability part. For the impermeability result, it is enough to set $\gamma=\min(\gamma^1,\dots,\gamma^M)$. Therefore, Theorem~\ref{main2} would remain unchanged.

In contrast, the case of inclusions distributed on the segment is more complicated. It is clear that we get the permeability if $ \frac{d_{\eps}}{\eps^{2+\frac1{\max\gamma^{p_{i,j}}}}} \to \infty$ and impermeability if $\frac{d_{\eps}}{\eps^{2+\frac1{\min \gamma^{p_{i,j}}}}} \to 0$, but it leaves aside many configurations. For example, if we distribute squares ($\gamma^1=\infty$) in the first half part of the segment and disks ($\gamma^2=1$) in the second part, then we guess that for $\eps^3 \ll d_{\eps} \ll \eps^2$ we obtain at the limit the Euler solution in the exterior of the half segment (i.e. the fluid cannot penetrate the part with the squares while it does not feel the presence of the disks).
\end{itemize}

Another possible extension (which is less easier) is the case where the obstacle has a corner i.e. $\gamma=0$. In this case, we could expect that the fluid is never perturbed by the porous medium if the obstacles are distributed on the segment. Actually, the cutoff function constructed in this paper should be sufficient in this setting, but the difficulty here is to perform the estimate of $w^{k,\eps}_{i,j}$. On the one hand, when $\gamma>0$, the boundary is $C^{1,\beta}$ which implies that the conformal mapping $\mathcal{T}$ from $\K^c$ to the exterior of the unit disk has a bounded derivative up to the boundary (Kellogg-Warschawski theorem). On the other hand, outside a corner, $D\mathcal{T}$ blows up in the neighborhood of the corner. The explicit form of this blow up is well-known by elliptic estimates in domains with piecewise smooth boundaries (see e.g. \cite{Kondratiev,Mazya}) and should be used to perform this extension.

Next, the same question could be investigated when the inclusions are distributed on the square.

\subsection{Infinite number of lines shrinking to the segment}

Another possible extension is the case where we consider several lines of obstacles, but in a thin region which shrinks to the segment. For example, let us introduce a parameter $\mu\in (0,1)$, and we consider $N_{2,\eps} := [(1/\eps)^\mu ]$ lines of obstacles separated by a vertical distance $\eps$ and by an horizontal distance $d_{\eps}$, i.e.
\begin{multline}\label{domain3}
 z_{i,j}^{\varepsilon}=(\tfrac\eps 2+(i-1)(\eps+d_{\varepsilon}),(j-1)2\eps)=(\tfrac\eps 2,0)+(\eps+d_{\eps})(i-1,0)+2\eps(0,j-1),\\ i=1,\dots, N_{\eps}, \ j=1,\dots, N_{2,\eps}.
\end{multline}
Therefore, the obstacles are included in the rectangle $[0,1]\times[-\eps/2,2 \eps^{1-\mu}]$ which shrinks to the unit segment as $\eps\to 0$.

In this situation, the only difference in the permeability analysis is that 
\[
 \sqrt{\sum_{i,j} 1} =\sqrt{N_{\eps}N_{2,\eps}} \leq \frac1{\eps^{\mu/2}}\min\Big(\sqrt{\frac2\eps},\sqrt{\frac2{d_\eps}}\Big),
\]
so using the same cutoff than in the case of a line, we get from \eqref{est:permeability} that
\begin{align*}
 \| K_{\R^2}[f]-v^\eps \|_{L^2(\Omega^\eps)} 
 &\leq C \|f\|_{L^1\cap L^\infty} \frac{ \sqrt{\eps}}{\eps^{\mu/2}} \Bigl(1+ \Big(\frac{\eps}{d_{\eps}}\Big)^{\frac{\gamma}{2(\gamma+1)}} \Big)\\
 &\leq C \|f\|_{L^1\cap L^\infty} \Big( \eps^{(1-\mu)/2} + \Big(\frac{\eps^{1+(1-\mu)(1+\frac1\gamma) }}{d_{\eps}}\Big)^{\frac{\gamma}{2(\gamma+1)}} \Big).
\end{align*}
This allows us to prove the permeability in the case where $\frac{d_{\eps}}{\eps^{1+(1-\mu)(1+\frac1\gamma) }} \to \infty$. For $\mu=0$, we recover the criterion for one line $\frac{d_{\eps}}{\eps^{2+\frac1\gamma }}$ whereas for $\mu=1$ we get the criterion corresponding to inclusions distributed on the square: $\frac{d_{\eps}}{\eps}$.

Following Section~\ref{imperm:square} with $N' = [(N_{2,\eps}-1)/2]$, it is also possible to adapt the impermeability result when $\frac{d_{\eps}}{\eps^{1+(1-\mu)(1+\frac1\gamma) }} \to 0$:
\begin{align*}
 (N's_{\eps})^{-1} \Big| \sum_{k=1}^{N'} \int_{-s_{\eps}}^{s_{\eps}} \int_{ \partial \mathcal{C}_{k}^\varepsilon(s)\setminus \partial \Omega^\varepsilon } \tilde p_{1,\eta} v^\eps \cdot n\, dxds \Big| 
 &\leq C\sqrt{2}\Big(\frac{d_{\eps}}{(N'+1)s_{\eps} \eps^2} \Big)^{1/2} \| p_{1,\eta} \|_{L^\infty(\Omega_2)} \| v^\eps \|_{L^2(\Omega_2^\varepsilon\cap \supp \tilde p_{1,\eta})}\\
 &\leq C\Big(\frac{d_\eps}{\eps^{1+(1-\mu)(1+\frac1\gamma)}}\Big)^{\frac\gamma{2(\gamma+1)}} \| p_{1,\eta} \|_{L^\infty(\Omega_2)} \| v^\eps \|_{L^2(\Omega_2^\varepsilon\cap \supp \tilde p_{1,\eta})}.
\end{align*}

Hence, we have the following.
\begin{corollary} Assume that $\K$ verifies {\rm (H1)-(H2)} and let $\mu\in [0,1)$.
Let $\omega_{0}\in L^\infty_{c}(\R^2)$ and $(u^\eps,\omega^\eps)$ be the global weak solution to the Euler equations \eqref{Euler1}-\eqref{Euler5} on 
$$\Omega^{\varepsilon}_{\mu}:=\R^2 \setminus \Big( \bigcup_{i=1}^{N_{\eps}}\bigcup_{j=1}^{N_{2,\eps}} \K_{i,j}^{\varepsilon}\Big)\quad \text{(with $ \K_{i,j}^{\varepsilon}$ defined in \eqref{domain1} and \eqref{domain3})},$$
with initial vorticity $\omega_{0}\vert_{\Omega^{\varepsilon}}$ and initial circulations $0$ around the inclusions (see \eqref{initial}).
\begin{enumerate}
 \item[(i)] If 
 $$
\frac{d_{\eps}}{\eps^{1+(1-\mu)(1+\frac1\gamma) }} \to \infty \quad \text{ for a sequence }\varepsilon \to 0
 $$
 then
 \begin{itemize}
\item $u^\eps \to u$ strongly in $L^2_{\loc}(\R^+\times\R^2)$ and $\omega^\eps \rightharpoonup {\omega}$ weak $*$ in $L^\infty(\R^+\times\R^2)$;
\item the limit pair $(u,{\omega})$ is the unique global solution to the Euler equations in the whole plane $\R^2$, with initial vorticity $\omega_{0}$.
\end{itemize}
 \item[(ii)] If 
 $$
\frac{d_{\eps}}{\eps^{1+(1-\mu)(1+\frac1\gamma) }} \to 0 \quad \text{ for a sequence }\varepsilon \to 0
 $$
 then there exists a subsequence such that
 \begin{itemize}
\item $u^{\eps} \rightharpoonup u$ weak $*$ in $L^\infty_{\loc}(\R^+; L^2_{\loc}(\R^2\setminus ([0,1]\times\{0\})))$ and $ \omega^{\eps} \rightharpoonup {\omega}$ weak $*$ in $L^\infty(\R^+\times(\R^2 \setminus ([0,1]\times\{0\})))$;
\item the limit pair $(u,{\omega})$ is a global weak solution to the Euler equations in $\R^2\setminus ([0,1]\times\{0\})$, with $u\cdot n=0$ on the boundary, with initial vorticity $\omega_{0}$ and initial circulation $0$ around the segment.
\end{itemize}
\end{enumerate}
\end{corollary}

\subsection{Non-uniform continuity of the Leray projection}

For the permeability part and the impermeability part, the $L^2$ space plays a special role, due to the orthogonality of the Leray projector in this framework. In this case, we state that the projector is uniformly continuous ($1$-lipschitz) no matter of the domain (even if there is more and more inclusions with tiny size). The natural question is to wonder if this projector is uniformly continuous in the $L^p$ framework for $p\neq 2$. Let us show that it is not the case for $p<2$ and $d_{\eps}=\varepsilon^\alpha$ with $\alpha>2$.

We consider the case of flat inclusions ($\gamma=\infty$) distributed on the segment. In the permeability part, the idea is to take advantage of the explicit formula of a correction $v^\varepsilon$ to prove that $K_{\R^2}[\omega^\varepsilon]-v^\varepsilon$ tends to zero in $L^2$ for $\alpha<2$ and to conclude because $u^\varepsilon-v^\varepsilon$ is the Leray projection in $\Omega^\varepsilon$ of $K_{\R^2}[\omega^\varepsilon]-v^\varepsilon$. However, for $p<2$ we can find some $\alpha>2$ for which $K_{\R^2}[\omega^\varepsilon]-v^\varepsilon$ tends to zero in $L^p$. This is an easy consequence of the analysis of Section~\ref{sec:permeability}: the convergence of $ \sum_{i} \varphi_{i,1}^\eps(x) (w_{i,1}^{3,\eps}+w_{i,1}^{4,\eps})$ to zero holds for any $p$, and considering a basic cutoff function (see e.g. \cite{BLM}) the second part could be estimated as
\begin{align*}
\Big\| \sum_{i} \nabla^\perp \varphi_{i,j}^\eps(x) (w_{i,j}^{1,\eps}+w_{i,j}^{2,\eps}) \Big\|_{L^2}\leq \frac{\eps}{d_{\eps}} \Big( \frac1{\eps} \eps d_{\eps}\Big)^{1/p}\leq \eps^{1-\alpha(1-\frac1p)}.
\end{align*}
We recover that for $p=2$, we should have $\alpha<2$, but for any $p\in [1,2),$ we can find $\alpha>2$ such that the right hand side term tends to zero. If the Leray projector is uniformly continuous in $L^p(\Omega^\varepsilon_{1})$, it would imply that the limit of $u^\varepsilon$ is the solution without influence of the porous medium. This is in contradiction with the point (ii) of Theorem~\ref{main1} where a wall appears.

Working on inclusions distributed on the square, we could also provide a counter example with $\alpha>1$. This is interesting because in the case of $\alpha=1$, the Leray projector is uniformly continuous in $L^p$ for any $p$ (see \cite{MasmoudiCRAS}).

\subsection{Paradox with the Kelvin theorem}

The Kelvin theorem states that the circulation around an inclusion is conserved for $t>0$. Therefore, if we consider an inclusion which shrinks to a point with zero initial circulation, at the limit the circulation around the point is zero (see e.g. \cite{ILL,LLL,Lopes}).

With zero initial circulation, the circulation of an inviscid flow around the unit segment (see \cite{Lac-IHP}) stays also zero. However, if we compute the curl of $u$, there is a measure supported on the curve $g(s,t)\delta$ where $\delta$ is the Dirac function on the segment $[0,1]\times\{0\}$. The density $g(s,t)$ is the jump of the tangential part and depends on time. The total circulation around the segment $\int g(\cdot,t)$ is equal to zero but $g$ is not zero in general. The presence of $g$ is important in order that $u=K_{\R^2}[\omega +g(t)\delta]$ is tangent to the boundary, while there is no reason that $K_{\R^2}[\omega]$ is tangent. Hence the density $g$ at a point $x\in [0,1]\times\{ 0 \}$ changes in order to counterbalance the normal part of $K_{\R^2}[\omega]$ and its total mass is zero. The connectedness of the segment allows this transfer of vorticity from one point to another of the segment.

If we think the segment as an infinite number of points, the Kelvin theorem can appear in contradiction with the part (ii) of Theorem~\ref{main1}. Nevertheless, the Kelvin theorem states only that the velocity at the point $x_{i,0}^\varepsilon \pm(\varepsilon,0)$ may be non zero (in particular if $u^{up}\neq u^{down}$) in order that the circulation around is zero. We have established that there is no average flux ($\int u^{up}\cdot e_{2} \varphi =\int u^{down}\cdot e_{2} \varphi = 0$ for any $\varphi \in C^\infty_{c}([0,1]\times \{0\})$), but there are maybe some vertical velocities close to the inclusions.

Actually, this remark raises another interesting question which is to understand the convergence with higher norms and to get the first order of $u^\varepsilon-u$.

\bigskip

\noindent
 {\bf Acknowledgements.} The first author is partially supported by the ANR Project DYFICOLTI grant ANR-13-BS01-0003-01 and by the Project ``Instabilities in Hydrodynamics'' funded by Paris city hall (program ``Emergences'') and the Fondation Sciences Math\'ematiques de Paris. The second author is partially supported by NSF grant DMS-1211806. 

The authors are grateful to the anonymous referees for their valuable comments on the first version of this article which led to a substantial improvement of this work.

\def\cprime{$'$}


\begin{thebibliography}{10}

\bibitem{Allaire90a}
G.~Allaire.
\newblock Homogenization of the {N}avier-{S}tokes equations in open sets
  perforated with tiny holes. {I}. {A}bstract framework, a volume distribution
  of holes.
\newblock {\em Arch. Rational Mech. Anal.}, 113(3):209--259, 1990.

\bibitem{Allaire90b}
G.~Allaire.
\newblock Homogenization of the {N}avier-{S}tokes equations in open sets
  perforated with tiny holes. {II}. {N}oncritical sizes of the holes for a
  volume distribution and a surface distribution of holes.
\newblock {\em Arch. Rational Mech. Anal.}, 113(3):261--298, 1990.

\bibitem{Allaire91}
G.~Allaire.
\newblock Homogenization of the {N}avier-{S}tokes equations with a slip
  boundary condition.
\newblock {\em Comm. Pure Appl. Math.}, 44(6):605--641, 1991.

\bibitem{AMBPT13}
A.~Bendali, M.~Fares, E.~Piot, and S.~Tordeux.
\newblock Mathematical justification of the {R}ayleigh conductivity model for
  perforated plates in acoustics.
\newblock {\em SIAM J. Appl. Math.}, 73(1):438--459, 2013.

\bibitem{BLM}
V.~Bonnaillie-No{\"e}l, C.~Lacave, and N.~Masmoudi.
\newblock Permeability through a perforated domain for the incompressible 2{D}
  {E}uler equations.
\newblock {\em Ann. Inst. H. Poincar\'e Anal. Non Lin\'eaire}, 32(1):159--182,
  2015.

\bibitem{Kondratiev}
M.~Borsuk and V.~Kondratiev.
\newblock {\em Elliptic boundary value problems of second order in piecewise
  smooth domains}, volume~69 of {\em North-Holland Mathematical Library}.
\newblock Elsevier Science B.V., Amsterdam, 2006.

\bibitem{CNS}
G.~Cardone, S.~A. Nazarov, and J.~Sokolowski.
\newblock Asymptotics of solutions of the {N}eumann problem in a domain with
  closely posed components of the boundary.
\newblock {\em Asymptot. Anal.}, 62(1-2):41--88, 2009.

\bibitem{CM82}
D.~Cioranescu and F.~Murat.
\newblock Un terme \'etrange venu d'ailleurs.
\newblock In {\em Nonlinear partial differential equations and their
  applications. {C}oll\`ege de {F}rance {S}eminar, {V}ol. {II} ({P}aris,
  1979/1980)}, volume~60 of {\em Res. Notes in Math.}, pages 98--138, 389--390.
  Pitman, Boston, Mass., 1982.

\bibitem{conca1}
C.~Conca.
\newblock \'{E}tude d'un fluide traversant une paroi perfor\'ee. {I}.
  {C}omportement limite pr\`es de la paroi.
\newblock {\em J. Math. Pures Appl. (9)}, 66(1):1--43, 1987.

\bibitem{conca2}
C.~Conca.
\newblock \'{E}tude d'un fluide traversant une paroi perfor\'ee. {II}.
  {C}omportement limite loin de la paroi.
\newblock {\em J. Math. Pures Appl. (9)}, 66(1):45--70, 1987.

\bibitem{Diaz99}
J.~I. D{\'\i}az.
\newblock Two problems in homogenization of porous media.
\newblock In {\em Proceedings of the Second International Seminar on Geometry,
  Continua and Microstructure (Getafe, 1998)}, volume~14, pages 141--155, 1999.

\bibitem{DAM12}
J.~Diaz-Alban and N.~Masmoudi.
\newblock Asymptotic analysis of acoustic waves in a porous medium: initial
  layers in time.
\newblock {\em Commun. Math. Sci.}, 10(1):239--265, 2012.

\bibitem{Galdi}
G.~P. Galdi.
\newblock {\em An introduction to the mathematical theory of the
  {N}avier-{S}tokes equations}.
\newblock Springer Monographs in Mathematics. Springer, New York, second
  edition, 2011.
\newblock Steady-state problems.

\bibitem{GV-Lac}
D.~G{\'e}rard-Varet and C.~Lacave.
\newblock The {T}wo-{D}imensional {E}uler {E}quations on {S}ingular {D}omains.
\newblock {\em Arch. Ration. Mech. Anal.}, 209(1):131--170, 2013.

\bibitem{GV-Lac2}
D.~G{\'e}rard-Varet and C.~Lacave.
\newblock The {T}wo-{D}imensional {E}uler {E}quations on {S}ingular {E}xterior
  {D}omains.
\newblock {\em To appear in Arch. Ration. Mech. Anal.}, 2015.

\bibitem{ILL}
D.~Iftimie, M.~C. Lopes~Filho, and H.~J. Nussenzveig~Lopes.
\newblock Two dimensional incompressible ideal flow around a small obstacle.
\newblock {\em Comm. Partial Differential Equations}, 28(1-2):349--379, 2003.

\bibitem{Kikuchi}
K.~Kikuchi.
\newblock Exterior problem for the two-dimensional {E}uler equation.
\newblock {\em J. Fac. Sci. Univ. Tokyo Sect. IA Math.}, 30(1):63--92, 1983.

\bibitem{Mazya}
V.~A. Kozlov, V.~G. Maz{\cprime}ya, and J.~Rossmann.
\newblock {\em Spectral problems associated with corner singularities of
  solutions to elliptic equations}, volume~85 of {\em Mathematical Surveys and
  Monographs}.
\newblock American Mathematical Society, Providence, RI, 2001.

\bibitem{Lac-IHP}
C.~Lacave.
\newblock Two dimensional incompressible ideal flow around a thin obstacle
  tending to a curve.
\newblock {\em Ann. Inst. H. Poincar\'e Anal. Non Lin\'eaire},
  26(4):1121--1148, 2009.

\bibitem{Lac-uni}
C.~Lacave.
\newblock Uniqueness for two-dimensional incompressible ideal flow on singular
  domains.
\newblock {\em SIAM J. Math. Anal.}, 47(2):1615--1664, 2015.

\bibitem{LLL}
C.~Lacave, M.~C. Lopes~Filho, and H.~J. Nussenzveig~Lopes.
\newblock Asymptotic behavior of 2d incompressible ideal flow around small
  disks.
\newblock In progress, 2015.

\bibitem{LaMa}
C.~Lacave and A.~Mazzucato.
\newblock The vanishing viscosity limit in the presence of a porous medium.
\newblock {\em arXiv preprint arXiv:1503.06554}, 2015.

\bibitem{PLL}
P.-L. Lions.
\newblock {\em Mathematical topics in fluid mechanics. {V}ol. 1}, volume~3 of
  {\em Oxford Lecture Series in Mathematics and its Applications}.
\newblock The Clarendon Press, Oxford University Press, New York, 1996.
\newblock Incompressible models, Oxford Science Publications.

\bibitem{LM99}
P.-L. Lions and N.~Masmoudi.
\newblock Une approche locale de la limite incompressible. (french) [a local
  approach to the incompressible limit].
\newblock {\em C. R. Acad. Sci. Paris S\'er. I Math.}, 329(5):387--392, 1999.

\bibitem{LionsMasmoudi}
P.-L. Lions and N.~Masmoudi.
\newblock Homogenization of the {E}uler system in a 2{D} porous medium.
\newblock {\em J. Math. Pures Appl. (9)}, 84(1):1--20, 2005.

\bibitem{Lopes}
M.~C. Lopes~Filho.
\newblock Vortex dynamics in a two-dimensional domain with holes and the small
  obstacle limit.
\newblock {\em SIAM J. Math. Anal.}, 39(2):422--436 (electronic), 2007.

\bibitem{Masmoudi02esaim}
N.~Masmoudi.
\newblock Homogenization of the compressible {N}avier-{S}tokes equations in a
  porous medium.
\newblock {\em ESAIM Control Optim. Calc. Var.}, 8:885--906 (electronic), 2002.
\newblock A tribute to J. L. Lions.

\bibitem{MasmoudiCRAS}
N.~Masmoudi.
\newblock Some uniform elliptic estimates in a porous medium.
\newblock {\em C. R. Math. Acad. Sci. Paris}, 339(12):849--854, 2004.

\bibitem{Mikelic91}
A.~Mikeli{\'c}.
\newblock Homogenization of nonstationary {N}avier-{S}tokes equations in a
  domain with a grained boundary.
\newblock {\em Ann. Mat. Pura Appl. (4)}, 158:167--179, 1991.

\bibitem{MikelicPaoli}
A.~Mikeli{\'c} and L.~Paoli.
\newblock Homogenization of the inviscid incompressible fluid flow through a
  {$2$}{D} porous medium.
\newblock {\em Proc. Amer. Math. Soc.}, 127(7):2019--2028, 1999.

\bibitem{MunnierRamdani}
A.~Munnier and K.~Ramdani.
\newblock Asymptotic analysis of a neumann problem in a domain with cusp.
  application to the collision problem of rigid bodies in a perfect fluid.
\newblock {\em arXiv preprint arXiv:1405.5446}, 2014.

\bibitem{Sanchez80}
E.~S{\'a}nchez-Palencia.
\newblock {\em Nonhomogeneous media and vibration theory}.
\newblock Springer-Verlag, Berlin, 1980.

\bibitem{SP1}
E.~S{\'a}nchez-Palencia.
\newblock Boundary value problems in domains containing perforated walls.
\newblock In {\em Nonlinear partial differential equations and their
  applications. {C}oll\`ege de {F}rance {S}eminar, {V}ol. {III} ({P}aris,
  1980/1981)}, volume~70 of {\em Res. Notes in Math.}, pages 309--325. Pitman,
  Boston, Mass., 1982.

\bibitem{Tartar80}
L.~Tartar.
\newblock Incompressible fluid flow in a porous medium: convergence of the
  homogenization process.
\newblock {\em {\it in} Nonhomogeneous media and vibration theory (E.
  S{\'a}nchez-Palencia)}, pages 368--377, 1980.

\bibitem{Wolibner}
W.~Wolibner.
\newblock Un theor\`eme sur l'existence du mouvement plan d'un fluide parfait,
  homog\`ene, incompressible, pendant un temps infiniment long.
\newblock {\em Math. Z.}, 37(1):698--726, 1933.

\bibitem{Yudo}
V.~I. Yudovi{\v{c}}.
\newblock Non-stationary flows of an ideal incompressible fluid.
\newblock {\em \u Z. Vy\v cisl. Mat. i Mat. Fiz.}, 3:1032--1066, 1963.

\end{thebibliography}
\end{document}